\documentclass[11pt,a4]{article}
\usepackage[colorlinks=true, pdfstartview=FitV, linkcolor=blue, citecolor=blue,
urlcolor=blue]{hyperref}
\usepackage[applemac]{inputenc}
\usepackage{amssymb,amsmath,amsthm,amsfonts,amscd,amstext}
\usepackage[mathscr]{eucal}
\usepackage{multirow}
\usepackage{xypic}
\usepackage{array}
\usepackage{graphicx}
\newtheorem{thm}{Theorem}[section]

\newtheorem{prop}[thm]{Proposition}

\newtheorem{df}[thm]{Definition\rm}
\newtheorem{dfthm}[thm]{Definition / Theorem}
\newtheorem{rem}{\it Remark\/}
\def\Aut{\operatorname {Aut}} 
\def\C{{\mathbb C}}  

\def\D{{\mathbb D}}
\def\deg  {{\rm deg}}
\def\der{\operatorname {Der}} 

\def\Det{\operatorname {Det}}

\def\dim{\operatorname{dim}}  

\def\E{{\mathbb E}}    

\def\Hom{\operatorname{Hom}} 
\def\id{\operatorname{id}}  

\def\ker{\operatorname{ker}}  

\def\OO{\mathcal O}  

\def\P{{\mathbb P}}  
\def\pr{\operatorname{pr}} 

\def\Q{{\mathbb Q}}   
\def\R{{\mathbb R}}    

\def\W{{\mathbb W}}    
\def\Y{{\mathbb Y}} 
\def\Z{{\mathbb Z}}    
\def\N{{\mathbb N}}    

\def\og{\leavevmode\raise.3ex\hbox{$\scriptscriptstyle\langle\!\langle$~}}

\def\fg{\leavevmode\raise.3ex\hbox{~$\!\scriptscriptstyle\,\rangle\!\rangle$}}
\begin{document}
\title{\bf Formal Groups, Witt vectors and Free Probability}
\author
{Roland Friedrich and John McKay}
\maketitle
\begin{abstract}
We establish a link between free probability theory and Witt vectors, via the theory of formal groups.  We derive an exponential isomorphism which expresses Voiculescu's free multiplicative convolution $\boxtimes$ as a function of the free additive convolution $\boxplus$. Subsequently we continue our previous discussion of the relation between complex cobordism and free probability. We show that the generic $n$th free cumulant corresponds to the cobordism class of the $(n-1)$-dimensional complex projective space. This permits us to relate several probability distributions from random matrix theory to known genera, and to build a dictionary. Finally, we discuss aspects of free probability and the asymptotic representation theory of the symmetric group from a conformal field theoretic perspective and show that every distribution with mean zero is embeddable into the Universal Grassmannian of Sato-Segal-Wilson.

\end{abstract}
\tableofcontents
\section{Introduction}
By considering the reduced free product of $C^*$-algebras, D.V. Voiculescu~\cite{V1985,VDN} was led to create free probability, a non-commutative probability theory where the central notion of independence is replaced by ``freeness". Although deep connections are known to exist between $C^*$-algebras, and more generally operator algebras and algebraic topology, these links have not (yet) pervaded (non-commutative) probability theory, as one might expect given the close links operator theory has with probability theory.  

In our previous article~\cite{FMcK2011}, we demonstrated among other things, a direct connection of one-dimensional free probability with multiplicative genera in complex cobordism via Voiculescu's $S$-transform~\cite{V1986}. The bridge Katsura, Shimizu and K. Ueno~\cite{KSU}  and J. Morava~\cite{Mo} had previously built, permits us to establish a link between free probability, conformal field theory (CFT) and the KP hierarchy.

Most notably, we observed a possible connection free probability has with the ring of Witt vectors~\cite{Hz} and, more generally, the theory of formal group laws.

First, in Theorem~\ref{LOG} our identification of the free additive and multiplicative convolution group with Witt vectors permits us derive a natural logarithm $\operatorname{LOG}$, which relates the two operations, analogously to the classical case. 

Second, the free additive and multiplicative convolution group define affine group schemes and formal group laws which are strongly isomorphic to those underlying Witt vectors. 

Then, for both the free additive and multiplicative convolution group we introduce their  representing Hopf algebras.

Third, the identification with Witt vectors 
permits us to define previously unknown operations on the space of distributions itself which are, a priori, not given by the underlying algebra, and which give rise to a commutative unital ring structure. In the case of actual probability measures positivity has to be taken into account, and then the ring structure reduces to a semiring, as shown in~\cite{F2019}.

From a different perspective the appearance of formal group laws in probability theory can be motivated as follows. 

A fundamental question in a general probability theory is, given two (non-commutative) random variables $a$ and $b$, how to express the moments, if possible at all, of $a+b$ and $ab$, solely from the moments of $a$ and $b$.

For free random variables $a,b$ in a non-commutative probability space, Voiculescu solved the problem of expressing the $n$th moment $m_n(a+b)$ and $m_n(ab)$ of $a+b$ and $ab$, respectively, as a function of the moments $m_i(a)$ and $m_j(b)$~\cite{V1985,V1986,V1987} by using Lie theoretic methods.

Once a notion of ``independence" is present, this can always be achieved on general grounds, and it naturally leads to formal group laws, cf. e.g.~\cite{FMcK2015}.

Related questions appear in algebraic topology~\cite{A,BMN} where independence corresponds geometrically to a product structure of spaces. So, 
for two complex line bundles $\xi_1,\xi_2$ over a space $X$, with $c_1$ the first Chern-Conner-Floyd class, the expression for $c_1(\xi_1\otimes\xi_2)$ is given by a formal group law over the ring $\Omega_U$, i.e. $F(x,y)\in \Omega_U^*[[x,y]]$, such that 
$$
c_1(\xi_1\otimes\xi_2)=F(c_1(\xi_1),c_1(\xi_2))=c_1(\xi_1)+c_1(\xi_2)+\sum_{m,n=1}^{\infty}\, a_{m,n}\, c_1(\xi_1)^m c_1(\xi_2)^n
$$
with $x:=c_1(\xi_1),y:=c_1(\xi_2)\in MU^2(X)$.

Since the first version of this article appeared as a preprint on the arXiv, things have considerably evolved, and the link between Witt vectors and free probability, turns out to be fruitful. 

First, let us mention the related but independent work of M. Mastnak and A. Nica~\cite{MN}, who pioneered  the use of Hopf algebraic methods in order to investigate the free multiplicative convolution group of Nica and Speicher~\cite{NS}, thereby also establishing the relation with the Hopf algebra of symmetric functions. 

By now, efficient and powerful Hopf and homotopy algebraic, operadic, and Lie theoretic methods, pursued by different groups, have pervaded all the different flavours of probability theory, in order to abstractly understand, e.g. the nature of the so-called ``moment-cumulant formul\ae". 

Let us mention the series of contributions from K. Ebrahimi-Fard and F. Patras, who started with~\cite{KP2015}, mainly from the perspective of shuffle and half-shuffle algebras, M. Schürmann and his students, e.g.~\cite{MSch2017}, and J.S. Park and his coauthors' homotopy probability theory, e.g.~\cite{DPT}. Finally, we generalised several of the results of the present article in~\cite{FMcK2013,FMcK2015,F2017,F2019}.

Second, the connection given in Theorem~\ref{LOG} 
between the free additive and multiplicative convolution was subsequently further explored  in the analytic category by several authors. 

G. Cébron~\cite{C2016} studied a semigroup version of our map EXP adapted to the analytic setting for freely infinitely divisible distributions and introduced a random matrix model of it. M. Anshelevich and O. Arizmendi~\cite{AA2016} considered  exponential homomorphisms for general non-commutative probabilities, with a focus on analytic applications. Finally, in~\cite{F2019} the analytic counterpart is established for semirings of freely infinitely divisible measures. 

The content of the rest of the paper is organised as follows. 

In Section 2, which is composed of five short subsections, we recall first basic definitions from free probability in a generalised form and then introduce the Lie group of automorphisms of the formal line, fixing the origin. Then we introduce basic definitions and results from the theory of (commutative) formal groups. Next we interpret Voiculescu's $R$- and $S$-transform as natural transformations and we derive our first Theorem~\ref{LOG}. The last part of this section is a reminder on free cumulants.  

In Section 3, which is subdivided into four parts, first we give the Hopf algebraic representation of the free additive and multiplicative convolution group. Then we relate these two abelian groups with Witt vectors and define a ring structure on non-commutative distributions. Finally, we discuss relations of (second order) free probability with Faber polynomials and the Grunsky matrix.

In Section 4, we describe relations free probability has with multiplicative genera. First we present the Fock space versions of Voiculescu and Haagerup of the $R$- and $S$-transform in a succinct form. Then we relate multiplicative genera with the notions previously introduced and provide examples of genera as probability measures. 

In the final Section 5, we relate free probability to conformal field theory (CFT) and integrable systems in both the Boson and Fermion picture, which leads to a powerful unifying perspective. 
 
\section{Free probability and formal group laws}
\subsection{Generalised free probability}
Let $k$ be a field of characteristic zero, $\mathbf{Alg}_k$ the category of associative, not necessarily unital,  $k$-algebras, and $\mathbf{cAlg}_k$ the associative and commutative unital $k$-algebras. For $A\in\mathbf{cAlg}_k$, we denote by $A^{\times}$ its group of units. 

First we recall some definitions from free probability~\cite{VDN} for general $k$-algebras, cf.~\cite{FMcK2013}.

\begin{df}
Let $\mathcal{A}\in\mathbf{Alg}_k$, with unit $1_{\mathcal{A}}$, and let $\phi:\mathcal{A}\rightarrow k$ be a pointed $k$-linear functional with $\phi(1_{\mathcal{A}})=1_k$. The pair $(\mathcal{A},\phi)$, is called a {\bf non-commutative $k$-probability space}, and elements $a\in\mathcal{A}$, {\bf random variables}.
\end{df}
\begin{df}
The {\bf law} or {\bf distribution} of a random variable $a\in\mathcal{A}$, is the unital $k$-linear functional $\mu_a\in k[[X]]$,
given by $a^0:=1_{\mathcal{A}}$, $\mu_a(X^n):=\phi(a^n)$, $n\in\N$.
The coefficients $(m_n(a))_{n\in\N}$ of the power series $\mu_a$ are the {\bf moments} of $a$.
\end{df}
If $\mathcal{A}$ is a $C^*$-algebra and the linear functional $\phi$ a {\bf state}, i.e.  continuous and positive $(\phi(a^*a)\geq0)$, then the pair $(\mathcal{A},\phi)$ is called a {\bf $C^*$-probability space}. Similarly, we have a $W^*$-probability space if $\mathcal{A}$ is a von~Neumann algebra and $\phi$ is weakly continuous.  

We denote the category of algebraic $k$-probability spaces by $\mathbf{Alg_kP}$.
\begin{df}[Freeness]
Let $(\mathcal{A},\phi)\in\mathbf{Alg_kP}$ and $I$ be an index set.  A family of sub-algebras $(\mathcal{A}_i)_{i\in I}$, with $1_{\mathcal{A}}\in \mathcal{A}_i\subset\mathcal{A}$, $i\in I$, is called {\bf freely independent} if for all $n\in\N^*$, and $a_k\in {\mathcal{A}_{i_k}\cap\ker(\phi)}$, $1\leq k\leq n$, 
$$
a_1\cdots a_n\in{\ker(\phi)},
$$
holds, given that $i_k\neq i_{k+1}$, $1\leq k< n$.
\end{df}
\begin{rem}
The requirement for freeness is that {\em consecutive indices} have to be distinct. Hence, $i_k=i_{k+2}$ is possible. Also note that the notion of freeness is purely algebraic.
\end{rem}

Voiculescu~\cite{V1985,V1986,V1987} showed that for a pair of free random variables $\{a,b\}$ in a non-commutative probability space, with distributions $\mu_a$ and $\mu_b$, respectively,  $\mu_{a+b}$ and $\mu_{ab}$ depend functionally only on $\mu_a$ and $\mu_b$. 

The {\bf additive free convolution} of $\mu_a$ and $\mu_b$, denoted as $\mu_a\boxplus\mu_b$, is equal to the distribution of $\mu_{a+b}$. Assign to $x_n$ and $y_n$ degree $n$. Then the operation $\boxplus$ is given by commutative, symmetric and homogeneous universal polynomials $(P_n(x_1,\dots,x_n,y_1,\dots,y_n))_{n\in\N^{*}}$, of degree $n$ with $\Z$-coefficients, such that
\begin{eqnarray}
\label{V_P_n}
\mu_{a}\boxplus\mu_b(X^n)&=&P_n(\mu_a(X),\dots,\mu_a(X^n),\mu_b(X),\dots,\mu_b(X^n))\\\nonumber
&=&\mu_a(X^n)+\mu_b(X^n)+\tilde{P}_n(\mu_a(X),\dots,\mu_a(X^{n-1}),\mu_b(X),\dots,\mu_b(X^{n-1})).
\end{eqnarray}

The {\bf multiplicative free
convolution} of $\mu_a$ and $\mu_b$, denoted by $\mu_a\boxtimes\mu_b$, is equal to the distribution of $\mu_{ab}$. Assign to $x_n$ and $y_n$ degree $n$. Then the operation $\boxtimes$ is given by commutative, symmetric and homogeneous universal polynomials $(Q_n(x_1,\dots,x_n,y_1,\dots,y_n))_{n\in\N^{*}}$, of bidegree $n$ with $\Z$-coefficients, such that for $n\geq2$,
\begin{eqnarray}
\label{V_Q_n}
\mu_{a}\boxtimes\mu_b(X^n)&=&Q_n(\mu_a(X),\dots,\mu_a(X^n),\mu_b(X),\dots,\mu_b(X^n))\\\nonumber
&=&\mu_a(X^n)\mu_b(X)^n+\mu_a(X)^n\mu_b(X^n)+\tilde{Q}_n(\mu_a(X),\dots,\mu_a(X^{n-1}),\mu_b(X),\dots,\mu_b(X^{n-1}))
\end{eqnarray}
and with $Q_1(x_1,y_1)=x_1y_1$.

As a consequence of the properties listed above, for $n=2$, we have 
$$
Q_2(x_1,x_2,y_1,y_2)=x_2 y_1^2+x_1^2 y_2-x_1^2 y_1^2.
$$ 
\subsection{The Lie group $\Aut(\OO)$}
In the approach to conformal field theory by Kawamoto, Namikawa, Tsuchiya and Yamada~\cite{KNTY} the spaces of formal power series, which we are going to introduce, are fundamental, cf. the discussion in the monograph~\cite{FB}. They have been also linked with complex cobordism by V. Bukhshtaber and A. Shukorov~\cite{BMN} and J. Morava~\cite{Mo} in connection with the Landweber-Novikov algebra.
In~\cite{FMcK2011} we added and discussed its subsequent appearance in renormalisation theory, cf.~\cite{FM}.  

For $A\in\mathbf{cAlg}_k$, we define functorially the following spaces:
\begin{eqnarray*}
\mathcal{O}_A:=A[[z]]& = &\left\{\sum_{n=0}^{\infty} a_n z^n\right\}:~\text{the ring of formal power series with coefficients in $A$.}  \\
\mathfrak{m}_A &: = & \left\{\sum_{n=1}^{\infty} a_n z^n\right\}:~\text{the (unique) maximal ideal of $\mathcal{O}_A$.} 
\end{eqnarray*}

The {\bf group of automorphisms} $\Aut({\OO}_A)$ of the $A$-algebra $\mathcal{O}_A\equiv A[[z]]$, or equivalently, the automorphisms of the formal line over $A$ which fix the origin,  can be described as the set of power series $a_1z+a_2z^2+\dots,$ with $a_1\in A^{\times}$. Its subgroup $\Aut_+({\OO}_A)$ corresponds to those automorphisms with $a_1=1$.

These spaces of automorphisms carry the structure of an infinite Lie group with corresponding  Lie algebras:
\begin{eqnarray*}
\Aut_+({\OO}_A) &  &\qquad\der_+({\OO}_A)=z^2A[[z]]\frac{d}{dz} \\
\cap\qquad & & \qquad\qquad\cap \\
\Aut({\OO}_A) & &\qquad\der_0({\OO}_A)=zA[[z]]\frac{d}{dz} \\
& & \qquad\qquad\cap \\
& & \qquad\der({{\OO}}_A)=A[[z]]\frac{d}{dz}
\end{eqnarray*}

By defining $\ell_n:=-z^{n+1}\frac{d}{dz}$, for $n\in\Z$, we obtain the commutation relations of the {\bf Witt algebra}
\begin{equation*}
\label{Witt_com}
[\ell_m,\ell_n]=(m-n)\ell_{m+n}~.
\end{equation*}
The non-negative part of the Witt algebra corresponds to $\der_0({\OO}_A)$ whereas  the 
Lie group $\Aut_+({\OO}_A)$ is represented by the  {\bf Faà di Bruno Hopf algebra},~cf. e.g.~\cite{FM}, which is a graded connected, commutative but not co-commutative Hopf algebra. 

\subsection{Formal group laws}
\label{FGL}
We recall from the theory of formal group laws several facts in order to discuss further the above results. References are, e.g.~\cite{A,BMN,Fro,Hz1978,Laz}. 

Let $R$ be a commutative ring with unit. An {\bf $n$-dimensional formal group law} over $R$ is an $n$-tuple of power series $F=(F_1,\dots,F_n)$ in $2n$ commuting variables $X=(x_1,\dots,x_n)$ and $Y=(y_1,\dots,y_n)$, such that for all $i=1,\dots,n$, $n\in\N^*$,
$$
F_i(X,Y)\in R[[x_1,\dots,x_n,y_1,\dots,y_n]],
$$
and the following axioms are satisfied:  
\begin{eqnarray}
\label{FG_neutral}
F(X,0)&=&X,\quad F(0,Y)=Y,\quad\text{(neutral element)},\\
\label{FG_associativity}
F(F(X,Y),Z)&=&F(X,F(Y,Z)),\quad\text{(associativity)}.
\end{eqnarray}
The formal group law is {\bf commutative} if 
\begin{equation}
\label{FG_commutativity}
F(X,Y)=F(Y,X)
\end{equation}
holds.
Relation~(\ref{FG_neutral}) implies that for every $i$, one has 
\begin{equation}
\label{polynomial_FG}
F_i(X,Y)=x_i+y_i+\sum_{|\mathbf{k}|,|\mathbf{l}|\geq1}c_{\mathbf{k},\mathbf{l}}(i)x_1^{k_1}\cdots x_n^{k_n} \cdot y_1^{l_1}\cdots y_n^{l_n},
\end{equation}
where $k_1+\cdots+k_n=|\mathbf{k}|$, $l_1+\cdots+l_n=|\mathbf{l}|$, $k_i,l_j\in\N$ and $c_{\mathbf{k},\mathbf{l}}(i)\in A$.

The axioms~(\ref{FG_neutral}) and~(\ref{FG_associativity}) imply the existence of a unique {\bf inverse} $\iota(X)$, i.e. an $n$-tuple of power series $\iota(X)=(\iota_1(X),\dots,\iota_n(X))$ with $\iota_i(X)\in R[[X]]$, for all $i$, such that
\begin{equation}
\label{FG_inverse}
F_i(X,\iota(X))=0=F_i(\iota(X),X).
\end{equation}

In the one-dimensional case, the inverse $\iota(x)$ of $F(x,y)$, is the power series $\iota(x)=\sum_{j=1}^{\infty} a_j x^j$ with $a_j\in A$, which satisfies relation~(\ref{FG_inverse}).

If we assign in an $N$-dimensional formal group law $F$ degree $j$ to $x_j$ and $y_j$, then we call it {\bf homogeneous} if the the power series $F_j$ in~(\ref{polynomial_FG}) is homogeneous of degree $j$. Hence, every homogeneous formal group law is of the form~(\ref{homogeneity}) and polynomial.

A {\bf homomorphism $f:F\rightarrow G$ of formal groups} $F$ and $G$ over $A\in\mathbf{cAlg}_k$, is a power series $f(x)\in xA[[x]]$ with $f(0)=0$ such that
$$
f(F(x,y))=G(f(x),f(y))~.
$$
It is an {\bf isomorphism} if $f(x)\in\Aut(\mathcal{O}_A)$, i.e. $f'(0)\in A^{\times}$, and it is a {\bf strict (strong) isomorphism} if $f(x)\in\Aut_+(\mathcal{O}_A)$, i.e. $f'(0)=1$.

The formal group law $F$ is {\bf linearisable} if there exists an $f(x)\in\Aut_+(\mathcal{O}_A)$, called the {\bf logarithm of the formal group law}, such that
\begin{equation*}
\label{log_form_group_law}
f(F(x, y)) = f(x)+f(y).
\end{equation*}

It can be expressed in terms of the {\bf invariant differential} $\omega(x)$ as 
\begin{equation}
\label{inv_diff}
df(x):=\left(\frac{\partial F(x,y)}{\partial y}\Big|_{y=0}\right)^{-1}dx,
\end{equation}
which defines an element of the cotangent space $\der^{\vee}_+({\OO}_A)$ to $\Aut_+({\OO_A})$.

The two basic formal group laws are the {\bf additive group law}:
$$
F_a(x,y):=x+y,
$$ 
and the {\bf multiplicative group law}:
$$
F_m(x,y):=x+y+xy.
$$  
For any $\Q$-algebra, $F_a$ and $F_m$ are strongly isomorphic, with $f(z)=\ln(1+z)$. 
In general, for a {\em commutative} $n$-dimensional formal group law we have
\begin{thm}[cf. e.g.~\cite{Fro,Laz}]
\label{Q-theorems}
Let $A\in\mathbf{cAlg}_{\Q}$ and $F=(F_1,\dots,F_n)$ be an $n$-dimensional {\em commutative} formal group law over $A$. 
Then $F$ is linearisable, i.e. it is (strongly) isomorphic to the $n$-dimensional {\bf additive group law} $F^n_a$, with
$$
(F^{n}_a)_i(x_1,\dots,x_n,y_1,\dots, y_n)=x_i+y_i,
\qquad\text{for $i\in\{1,\dots,n\}$}.
$$
\end{thm}
\subsection{The $R$- and $S$-transform}
For $A\in\mathbf{cAlg}_k$, we define functorially the following sets of $A$-valued distributions: 
\begin{eqnarray*}
\Sigma^{}(A) & :=&\{(a_1, a_2,a_3,\dots), a_i\in A \}, \\
\Sigma^{\times}(A) & := &\{(a_1, a_2,a_3,\dots), a_1\in A^{\times} \}, \\
\Sigma_0(A)&:=&\{(0,a_2,a_3,\dots)\},\\
\Sigma_1(A)&:=&\{(1, a_2,a_3,\dots)\}, 
\end{eqnarray*}
which satisfy the strict inclusions:
$\Sigma_1(A)\subset\Sigma^{\times}(A)\subset\Sigma(A)$.

In [\cite{VDN} Theorems 3.2.3.and 3.6.3.] Voiculescu establishes the existence of two infinite complex Lie groups, induced by the free additive and multiplicative convolution.  This can be generalised and will be further discussed in~Section~\ref{Hopf algebras}.

\begin{thm}
\label{Alg_Groups}
 The pairs $(\Sigma,\boxplus)$ and $(\Sigma^{\times},\boxtimes)$ define pro-affine commutative group schemes over $k$, with the group laws $\boxplus$ and $\boxtimes$ given by the universal polynomials  $P_{\boxplus}:=(P_n)_{n\in\N^*}$ in~(\ref{V_P_n}) and  $Q_{\boxtimes}:=(Q_n)_{n\in\N^*}$ in~(\ref{V_Q_n}), respectively.
\end{thm}

For every $A\in\mathbf{cAlg}_k$ the group $\Sigma^{\times}(A)$ is the semi-direct product of the maximal torus $A^{\times}$ and the connected subgroup $\Sigma_1(A)$, viz. 
$$
\Sigma^{\times}(A)=A^{\times}\ltimes\Sigma_1(A).
$$

\begin{prop}
For every $A\in\mathbf{cAlg}_k$ and (fixed) $N\in\N^*$, the following holds: 
\begin{enumerate}
\item the polynomials $(P_n)_{n\in\N^*}$ in~(\ref{V_P_n}) define a commutative $N$-dimensional formal group law $F^N_{\boxplus}:=(P_1,\dots, P_N)$ over $\Z$, such that $P_{j}(X,Y)\in \Z[x_1,\dots,x_j,y_1,\dots,y_j]$ and
\begin{equation}
\label{homogeneity}
P_j(X,Y)=x_j+y_j+\sum_{|\mathbf{k}|,|\mathbf{l}]\geq1}c_{\mathbf{k},\mathbf{l}}(j)x_1^{k_1}\cdots x_{j-1}^{k_{j-1}}y_1^{l_1}\cdots y_{j-1}^{l_{j-1}},
\end{equation}
which satisfy the homogeneity condition $(k_1+l_1)+2(k_2+l_2)+\cdots+(j-1)(k_{j-1}+l_{j-1})=j$.
\item 
For every $A\in\mathbf{cAlg}_{\Q}$, the formal group law $F_{\boxplus}^N$ over $A$ is strongly isomorphic to the $N$-dimensional additive group law $F^N_a$.
\end{enumerate}
\end{prop}

The actual construction of the logarithm for the free additive convolution was done by Voiculescu~[\cite{V1985}, \cite{V1986}  and \cite{VDN} Section~3.2], and it is called the $R$-transform. We (re)state its content in intrinsic terms first.  Let $(\mathbb{A}^N,+)$ denote the $N$-dimensional additive group scheme.

\begin{dfthm}[Voiculescu's $R$-transform]
The $R$-transform is a natural isomorphism, i.e. for every $A\in\mathbf{cAlg}_k$  there exists an isomorphism of abelian groups
$$
R_A:(\Sigma(A),\boxplus)\rightarrow (A^{\N^*},+),
$$
such that for all $f,g\in\Sigma(A)$
$$
R_{A}(f\boxplus g)=R_A(f)+R_A(g),
$$
holds.
\end{dfthm}
Its definition is as follows, cf.~~[\cite{VDN}~p. 23-24].
For an $a\in\Sigma(A)$, with $a=(a_1,a_2,a_3,\dots)$ and $a_0:=1$, its {\bf Cauchy transform} is the formal power series 
\begin{equation*}
G_a(z):=\sum_{n=0}^{\infty} a_n\, z^{-(n+1)}.
\end{equation*}

The compositional inverse of $G_a(z)$, is the Laurent series 
\begin{equation*}
\label{invCauchy}
G_a^{-1}(z)=\frac{1}{z}+\sum_{n=0}^{\infty}\alpha_n z^n.
\end{equation*}
Then the {\bf ${\mathcal{R}}$-transform} of $a$ is the formal power series
\begin{equation*}
\label{R_Voiculescu}
\mathcal{R}_{a}(z):=G_a^{-1}(z)-\frac{1}{z}=\sum_{n=0}^{\infty}\alpha_n z^n. 
\end{equation*}
Often, it is more convenient to consider instead the series
\begin{equation*}
\label{R-trafo}
R_a(z):=z\mathcal{R}_a(z), 
\end{equation*}
where we use a simplified notation, i.e. we do not indicate the dependence on the $k$-algebra $A$. Further, depending on the context, we shall perceive the $R$-transform either as an $\mathbb{A}^{\N^*}$-valued natural transformation or as a formal power-series valued functor.   

\begin{rem}
If one expresses initial coefficients $\{b_k\}$ in the
{\bf Lagrange inversion} of series of the form $f(q) = 1/q+\sum_{k=0}^{\infty} a_k q^k$ into the reverted series $q = \sum_{k=0}^{\infty} b_k/f^{k+1}$, then each $b_k$ is a polynomial
over $\N$ in the original coefficients. This follows from an
interpretation of the products in the isobaric $b_k$ in
terms of {\bf Dyck paths}. Further, if the original series
has coefficients that are values of a {\bf character} on some
group element, $g$, then so are the $b_k$-polynomials,~\cite{FrdMcK}.
\end{rem}

For $A\in\mathbf{cAlg}_k$, the {\bf units} of the ring $A[[z]]$ are given by
$$
A[[z]]^{\times}=A^{\times}(1+zA[[z]]).
$$
Let $\bullet[[z]]^{\times}$ denote the functor $A\mapsto A[[z]]^{\times}$ and $\Lambda(A)\subset A[[z]]^{\times}$ the subset $\Lambda(A):=1+zA[[z]]$.

To every $f\in\Sigma^{\times}(A)$, with $f=(a_1, a_2,a_3,\dots)$ corresponds a power series $f(z)\in\Aut({\OO}_A)$ of the form
$$
f(z)=a_1z+a_2z^2+a_3 z^3+\dots,\qquad a_1\in A^{\times},
$$
with compositional inverse $f^{-1}(z)\in\Aut({\OO}_A)$.

\begin{dfthm}[Voiculescu's $S$-transform]
The {\bf $S$-transform} is, for every $A\in\mathbf{cAlg}_k$ and $f\in\Sigma^{\times}(A)$,  a natural isomorphism of abelian groups, 
$$
S_A:(\Sigma^{\times}(A),\boxtimes)\rightarrow (A[[z]]^{\times},\cdot)
$$ 
which is given by 
\begin{equation*}
\label{S_trafo}
S_A(f):=\frac{1+z}{z}f^{-1}(z).
\end{equation*}
For $f,g\in\Sigma^{\times}(A)$ it satisfies
$$
S_A(f\boxtimes g)=S_A(f)\cdot S_A(g).
$$
\end{dfthm}
In concrete terms the $S$-transform can be calculated as follows, cf.~[\cite{VDN} Theorem and Definition 3.6.3]: for $a\in\Sigma^{\times}(A)$, with $a=(a_1,a_2,a_3,\dots)$, $a_1\in A^{\times}$ the associated {\bf moment series} or {\bf $\mathcal{M}$-transform} is
\begin{equation}
\label{moment_series}
\mathcal{M}_a(z):=\sum_{n=0}^{\infty} a_{n+1}\, z^n,
\end{equation}
and $M_a(z):=z\mathcal{M}_a(z)\in\Aut({\OO}_A)$, with compositional inverse $M^{-1}_a(z)$. The {\bf $S$-transform} of $a$, is the formal power series
\begin{equation}
\label{S-moment-trafo}
S_{a}(z):=\frac{1+z}{z}M_{a}^{-1}(z), 
\end{equation}
where again we do not indicate the dependence on the $k$-algebra $A$.

Let us now compare the $R$- and $S$-transform with the integral transforms in classical probability theory. 
For of a positive real random variable $X$, i.e. $X\geq0$, the {\bf Mellin transform} $\mathfrak{M}$, cf. e.g.~\cite{GS}, is defined as
$$
\mathfrak{M}_X(t):=\E[X^{it}]=\int_0^{\infty} x^{it}d\mu_X(x),
$$
with $t\in\R$ and the convention $0^{it}:=0$, for all $t$.

Then, the following relations hold for $X,Y$ classically independent and $a,b$ free non-commutative random variables, respectively:
\[
\begin{tabular}{l|l|l} $X,Y$~{\bf independent}& convolution & transform (Fourier / Mellin)  \\\hline additive & $\mu_{X+Y}=\mu_X\ast\mu_Y$ & $\mathcal{F}_{X+Y}(t)=\mathcal{F}_X(t)\cdot\mathcal{F}_Y(t)$  \\multiplicative & $\mu_{X\cdot Y}=\mu_X\ast_{\operatorname{m}}\mu_Y$ & $\mathfrak{M}_{X\cdot Y}(t)=\mathfrak{M}_X(t)\cdot\mathfrak{M}_Y(t)$   \\\hline\hline 
$a,b$~{\bf free} & convolution & transform ($R$ / $S$)  \\\hline additive & $\mu_{a+b}=\mu_a\boxplus\mu_b$ & $R_{a+b}(z)=R_a(z)+R_b(z)$  \\multiplicative & $\mu_{a\cdot b}=\mu_a\boxtimes\mu_b$ & $S_{a\cdot b}(z)=S_a(z)\cdot S_b(z)$ \end{tabular}
\]
where we used the notation $\ast$ and $\ast_m$ for the respective convolution operations and $R_a(z):=R_{\mu_a}(z)$ and $S_a(z):=S_{\mu_a}(z)$.

Furthermore, by taking the natural logarithm of the Fourier, Mellin and $S$-transform, which are maps from the convolution group to the (co)-tangent space, the corresponding operations are linearised. 

In fact, in the free case there exists also a functional relation between $\boxtimes$ and $\boxplus$, analogous as for classically independent random variables between $\cdot$ and $+$.
\begin{thm}
\label{LOG}
Let $k$ be a field of characteristic zero and $A\in\mathbf{cAlg}_k$. There exists a natural group isomorphism, $\operatorname{LOG}:(\Sigma_1, \boxtimes)\rightarrow(\Sigma_0,\boxplus)$, such that for $f,g\in\Sigma_1(A)$ 
\begin{equation*}
\label{formula_free_mult_conv}
\operatorname{LOG}_A(f\boxtimes g)=\operatorname{LOG}_A(f)\boxplus\operatorname{LOG}_A(g),
\end{equation*}
holds. 
\end{thm}
\begin{proof}
One possibility to proof the statement is by composing the group isomorphisms in the diagram:
\[
\begin{xy}
  \xymatrix{
      (\Lambda(A), \cdot) \ar[d]^{z\frac{d}{dz}\ln} & (\Sigma_1(A), \boxtimes) \ar[l]_{S_A}\ar[d]^{=:\operatorname{LOG}_A} \\
                              (zA[[z]], +)\ar[r]^{\mathcal{R}_A^{-1}} & (\Sigma_0(A),\boxplus) 
               }
\end{xy}
\]

\end{proof}

\begin{rem}
If we are dealing with actual measures positivity matters, and then one obtains semi-groups instead of groups, cf.~\cite{F2019}. 
\end{rem}

\subsection{Free cumulants}

Let us recall the groundbreaking {\bf combinatorial interpretation} of the $R$- and $S$-transform given by R.~Speicher~\cite{S} and then jointly further developed with A.~Nica~\cite{NS}.

Let $S$ be a finite set, e.g. $\{1,\dots,n\}$. A {\bf partition} $\pi$ of $S$ is a collection $\pi=\{V_1,\dots, V_r\}$ of non-empty, pairwise disjoint subsets $V_i\subset S$, called the {\bf blocks} of $\pi$, whose union is again $S$, i.e. $S=\coprod_{i=1}^r V_i$. Let $\mathcal{P}(S)$ denote the set of all partitions of $S$, and  write $\mathcal{P}(n)$ if $S=\{1,\dots,n\}$. 
On $\{1,\dots,n\}$ there exists an equivalence relation $\sim_{\pi}$: namely, let $1\leq p,q\leq n$, and set
$$
p\sim_{\pi} q:\quad\Leftrightarrow\exists i: p,q\in V_i, \quad\text{i.e. $p,q$ are in the same block $V_i$ of $\pi$.}
$$
A partition $\pi$ of $S=\{1,\dots,n\}$ is called {\bf crossing} if there exist $p_1<q_1<p_2<q_2$ in $S$ such that $p_1\sim_{\pi} p_2$, $q_1\sim_{\pi} q_2$ but $p_2\nsim_{\pi} q_1$. A {\bf non-crossing} partition, as introduced by {\bf Kreweras}~\cite{S}, is a partition $\pi$ which is not crossing. The set of all non-crossing partitions of $\{1,\dots,n\}$ is denoted by $\operatorname{NC}(n)$.
These notions can be generalised to any finite, {\em totally ordered} set. 

For $\mathcal{A}\in\mathbf{Alg}_k$, let $(\mathcal{A},\phi)$ be a non-commutative probability space. The {\bf free cumulants} $\kappa_n$, are $k$-multilinear functionals $\kappa^n:\mathcal{A}^n\rightarrow k$ $(n\in\N)$ which are inductively defined by the {\bf moment-cumulant formula}:
\begin{equation*}
\kappa_1(a):=\phi(a),\qquad \phi(a_1\cdots a_n)=\sum_{\pi\in\operatorname{NC}(n)}\kappa_{\pi}[a_1,\dots,a_n]~,
\end{equation*}
where for $\pi=\{V_1,\dots, V_r\}\in\operatorname{NC}(n)$
$$
\kappa_{\pi}[a_1,\dots, a_n]:=\prod_{j=1}^r \kappa_{V_j}[a_1,\dots, a_n]~,
$$
and for $V=(v_1,\dots, v_s)$
$$
\kappa_V[a_1,\dots, a_n]:=\kappa_s[a_{v_1},\dots,a_{v_s}]~.
$$

For $a\in\mathcal{A}$ and $n\in\N^{*}$, one sets $\kappa_n(a):=\kappa_n(a,\dots,a)$ and calls the resulting series $(\kappa_n(a))_{n\in\N^*}$ the series of {\bf free cumulants of $a$}.

\begin{rem}
The above construction continues to hold for algebraic probability spaces with values in $\mathbf{cAlg}_k$, however, the original notion of freeness has to be appropriately redefined, as a generic random variable can not be centred,  cf.~\cite{FMcK2013,FMcK2015}.

Free cumulants are intrinsically Lie theoretic quantities~\cite{V1985}, which is also true for cumulants in other non-commutative probability theories,  cf.~\cite{KP2015,FMcK2015}.
\end{rem}

For $m\in\Sigma(A)$, with $m=(m_1,m_2,m_3,\dots)$ the free cumulants, $\kappa(m)=(\kappa_1(m),\kappa_2(m),\kappa_3(m),\dots)$ are related according to [\cite{NS} p. 209, equation~(12.2)] by the algebro-combinatorial formula
\begin{equation*}
\label{moment-cumulant_formula}
m_n=\sum_{\pi\in \operatorname{NC}(n)}\kappa_{\pi}(m).
\end{equation*}

For $a\in\Sigma(A)$, the coefficients of its $R$-transform are given by the free cumulants as
\begin{equation}
\label{R-cumul}
R_{a}(z)=\sum_{n=1}^{\infty} \kappa_n(a) z^n.
\end{equation}
If $k_1(a)\neq0$, the compositional inverse $R_{a}^{-1}(z)$ of~(\ref{R-cumul}) exists as an element of $\Aut(\OO_A)$.

Following~[\cite{NS} Definition 18.15.], the $S$-transform of an $a$ with $\kappa_1(a)\neq0$, can be expressed in terms of the $R$-transform as
\begin{equation}
\label{RScumul}
S_{a}(z)=\frac{R^{-1}_{a}(z)}{z},
\end{equation}
which according to~[\cite{NS} Formula(18.11)] is called the {\bf $\mathcal{F}$-transform}.
\section{Witt vectors and Hopf algebras}
\subsection{The affine group schemes}
\label{Hopf algebras}
We use two basic observations, namely commutative Hopf algebras are dual to affine group schemes ({\bf Cartier duality}) and an $n$-dimensional formal group law over $k$ defines a smooth formal $k$-group, and vice versa. References are~\cite{D1972,Hz1978,W1979}

In detail, for $k$ a commutative unital ring, an {\bf affine group scheme over $k$} or a {\bf $k$-group (functor)}, is a covariant functor $G:\mathbf{cAlg}_k\rightarrow\mathbf{Grps}$ which is representable, i.e. there exists a commutative but not necessarily co-commutative Hopf algebra $H$ over $k$, such that for all $A\in\mathbf{cAlg}_k$,
$$
G(A)\simeq\Hom_{\mathbf{cAlg}_k}(H,A),\quad G=\operatorname{Spf}(H)
$$
holds, as convolution groups.

Let us now assume that $k$ is a field of characteristic zero.
The {\bf coordinate algebra} $H=k[G]$ (or $\mathcal{O}(G)$), is {\bf co-connected} if there exists a chain $(C_n)_{n\in\N}$ of subspaces $C_n\subset H$ (filtration), such that $C_0=k$, $C_0\subseteq C_1\subseteq C_2\subseteq\dots$, $\bigcup_{n\in\N}C_n=H$ and
$$
\Delta(C_n)\subset\sum^n_{i=0}C_i\otimes C_{n-i}.
$$
An affine group scheme $G$ over $k$ is {\bf (pro)-unipotent} if its representing Hopf algebra $k[G]$ is co-connected, cf.~[\cite{W1979} Theorem p.64].

The formal scheme $G=\operatorname{Spf}(H)$ is {\bf connected} if $H$ is a local ring and a connected formal $k$-group $G=\operatorname{Spf}(H)$ is {\bf smooth} if $H$ is isomorphic to a {\bf power-series algebra} $k[[x_1,\dots,x_n]]$ in $n$ commuting variables. Every smooth formal $k$-group is connected, because $k[[x_1,\dots,x_n]]$ is a local ring. 

An $n$-dimensional formal group law $F(X,Y)$ over a $k$-algebra $A$ defines a co-product by:
\begin{eqnarray*}
\Delta_F:A[[x_1,\dots,x_n]]&\rightarrow&A[[x_1,\dots,x_n]]\hat{\otimes}A[[x_1,\dots,x_n]]\\
x_i&\mapsto&F_i(x_1\otimes 1,\dots,x_n\otimes1,1\otimes x_1,\dots,1\otimes x_n)=\\
&&x_i\otimes 1+1\otimes x_i+\text{(terms of degree $\geq3$)}.
\end{eqnarray*}
Together with the co-identity (augmentation) 
$$
\epsilon:A[[x_1,\dots,x_n]]\rightarrow A,\quad x_i\mapsto 0,
$$
and the co-inverse (antipode) $\alpha$, given by $X\mapsto\iota(X)$, where $\iota(X)$ is the inverse~(\ref{FG_inverse}), the $A$-algebra $A[[x_1,\dots,x_n]]$ becomes a commutative but not necessarily co-commutative Hopf algebra. 
For a smooth formal group $G=\operatorname{Spf}(H)$, $H=k[[x_1,\dots,x_n]]$, the co-multiplication $\Delta:H\rightarrow H\hat{\otimes}H$ is determined by the values 
$$
\Delta(x_i)=F_i(x_1\otimes1,\dots,x_n\otimes1,1\otimes x_1,\dots,1\otimes x_n),
$$
and hence by an $n$-tuple of power series $F=(F_1,\dots,F_n)$ in $2n$ variables. 

For $X_i:=x_1\otimes1$, $Y_i:=1\otimes x_i$, $X=(X_1,\dots,X_n)$ and $Y=(Y_1,\dots,Y_n)$, we have 
$$
k[[x_1,\dots,x_n]]\hat{\otimes}k[[x_1,\dots,x_n]]\cong k[[X,Y]].
$$
Then $F(X,Y)\in k[[X,Y]]$ and it has to satisfy the axioms of a formal group law ~(\ref{FG_neutral}) (neutral element), (\ref{FG_associativity}) (associativity) but not necessarily (\ref{FG_commutativity}) (commutativity).

The free polynomial algebra $\Z[x_1,x_2,x_3,\dots]$, in countably many commuting variables $x_i$, becomes a {\bf graded connected} algebra if we define $\operatorname{deg}(x_i):=i$, $i\in\N^*$, with the direct sum decomposition  
$$
\Z[x_1,x_2,x_3,\dots]=\bigoplus_{n=0}^{\infty} H_n,
$$
into homogeneous subspaces $H_j$, where $H_0= \Z\cdot 1_{\Z}$ and for $n\geq 1$, $H_n$  is the $\Z$-linear span of all monomials $x_{i_1}\cdots x_{i_m}$ of degree $i_1+\cdots+ i_m=n$, satisfying $H_i\cdot H_j\subseteq H_{i+j}$.  

We obtain a Hopf algebra $H_{\boxplus}$, by defining the co-unit $\varepsilon_{\boxplus}$ as $\varepsilon_{\boxtimes}(x_i):=0$ for $i\in\N^*$, the co-multiplication $\Delta_{\boxplus}$ as $\Delta_{\boxplus}(1)=1\otimes 1$, and on generators $x_n$, as
\begin{equation}
\Delta_{\boxplus}(x_n):=P_n(x_1\otimes1,\dots, x_{n}\otimes1,1\otimes x_1,\dots,1\otimes x_{n})\in\bigoplus_{j=0}^n H_j\otimes H_{n-j}
\end{equation}
with the universal polynomials $(P_n)_{n\in\N^*}$ from~(\ref{V_P_n}), and where the last inclusion follows from expression~(\ref{homogeneity}). The antipode $\alpha_{\boxplus}:H_{\boxplus}\rightarrow H_{\boxplus}$ is given by the inverse $\iota_{\boxplus}$ of the formal group law $F_{\boxtimes}$, and is recursively determined, starting with $\alpha_{\boxplus}(x_1)=\iota_{\boxplus}(X)=-x_1$. 

The following holds in general:
\begin{prop}
Let $k$ be a field of characteristic $0$, and $F$ an $N$-dimensional, $N\in\N^*\cup\{+\infty\}$, homogeneous formal group law over $k$. Then the associated Hopf algebra is (pro)-unipotent, and the co-product satisfies
\begin{equation*}
\Delta_F(x_i)=x_i\otimes 1+1\otimes x_i+\sum_j p_{ij}\otimes q_{ij},
\end{equation*}
for $i=1,\dots,N$, and polynomials $p_{ij},q_{ij}\in k[x_1,\dots, x_{i-1}]$.
\end{prop}

We note that for any field $k$ of characteristic zero, one has $\Z\subset k$.
\begin{prop}[Additive case]
The Hopf algebra $(k[x_1,x_2,x_3,\dots], \varepsilon_{\boxplus}, \Delta_{\boxplus}, \alpha_{\boxplus})$, induced by the free additive convolution $\boxplus$, is graded connected, co-connected and co-commutative. It represents the pro-unipotent formal group scheme $(\Sigma,\boxplus)$, and its $\C$-points correspond to Voiculescu's original infinite Lie group $(\Sigma(\C),\boxplus)$. 
\end{prop}

Let $k[x^{-1},x]]$ denote the ring of formal {\bf Laurent series}. 
\begin{prop}[Multiplicative case]
The pro-affine $k$-group scheme $(\Sigma^{\times},\boxtimes)$ is represented by 
the Hopf algebra $k[x_1^{-1},x_1,x_2,x_3,\dots]$, with the co-commutative co-product given by
$$
\Delta_{\boxtimes}(x_n):=x_n\otimes x_1^n+x_1^n\otimes x_n+\tilde{Q}_n(x_1\otimes 1,\dots,x_{n-1}\otimes 1,1\otimes x_1,\dots,1\otimes x_{n-1}),
$$
for $n\in\N^*$ with $n\geq2$, the polynomial ${\tilde{Q}}_n$ from~(\ref{V_Q_n}), and $x_1^{\pm1}$ being group-like. The co-unit satisfies $\varepsilon_{\boxtimes}(x_1^{\pm1})=\pm1_k$ and $\varepsilon_{\boxtimes}(x_n)=0$, for $n\geq2$. 
The antipode $\alpha_{\boxtimes}$ is calculated recursively, starting with $\alpha_{\boxtimes}(x_1^{\pm1})=x_1^{\mp1}$, and from the (general) identity $(\alpha_{\boxtimes}(x_n))(\bullet) = x_n(\iota_{\boxtimes}(\bullet))$, where $\iota_{\boxtimes}(\bullet)$ is the group inverse of a group element $\bullet$, and $x_n$ the $n$th coordinate function.
\end{prop}

Let $\mathbb{A}^N$, $N\in\N^*\cup\{+\infty\}$, denote the {\bf (pro) affine $N$-space over $k$}, i.e. the functor $\mathbb{A}^N:\mathbf{cAlg}_k\rightarrow\mathbf{Sets}$, with $A\mapsto A^N$.

The natural isomorphism 
\begin{eqnarray}
\label{translation}
\varphi_A:\mathbb{A}^{N}(A)&\rightarrow&\Sigma_1^{N+1}(A),\\\nonumber
a&\mapsto&(1,a_1+1,a_2+1,a_3+1,\dots,a_N+1),
\end{eqnarray}
for $A\in\mathbf{cAlg}_k$, induces, via the pullback, an universal commutative polynomial formal group law $F_{\tilde{\boxtimes}}:=\varphi^*(\boxtimes)$ over $k$, such that $(\mathbb{A}^{N},F_{\tilde{\boxtimes}})$ defines a $k$-group functor for all $N\in\N^*\cup\{+\infty\}$ and $1\leq n\leq N$, with
$$
F_{\tilde{\boxtimes},n}(x_1,\dots,x_n, y_1,\dots, y_n)=Q_{n+1}(1,x_1+1\dots,x_n+1, 1,y_1+1,\dots, y_n+1)-1.
$$
If we set $\deg(x_i)=\deg(y_i):=i$, then $F_{\tilde{\boxtimes},n}$ is homogeneous of degree $n$, which follows from the relation non-crossing partitions have with their Kreweras complement, cf.~\cite{FMcK2015,MN}.
Further, according to Theorem~\ref{Q-theorems}, the logarithm
\begin{equation}
\label{log_F_plus}
\log_{F_{\tilde{\boxtimes}}}:(\mathbb{A}^{N},F_{\tilde{\boxtimes}})\rightarrow(\mathbb{A}^{N},+)
\end{equation}
exists. 
\begin{prop}[Linearisation of $\boxtimes$] 
\label{main-iso}
There exists a natural isomorphism of $k$-group functors
$$
\log_{\boxtimes}:(\Sigma_1,\boxtimes)\rightarrow(\mathbb{A}^{\N^*},+),
$$ 
which for $A\in\mathbf{cAlg}_k$, is given by the diagram:
\[
\begin{xy}
  \xymatrix{
 (\Sigma_1(A),\boxtimes)\ar[r]^{{\varphi_A^{-1}}}\ar[dr]_{\log_{\boxtimes}:=\log_{F_{\tilde{\boxtimes}}}\circ\varphi_A^{-1}\qquad}    & (\mathbb{A}^{\N^*}(A),F_{\tilde{\boxtimes}})\ar[d]^{\log_{F_{\tilde{\boxtimes}}}}\\
          & (\mathbb{A}^{\N^*}(A),+)
               }
\end{xy}
\]
with $\varphi^{-1}_A$ the inverse morphism of~(\ref{translation}) and $\log_{F_{\tilde{\boxtimes}}}$ in~(\ref{log_F_plus}).
\end{prop}

In summary, for $(\mathcal{A},\phi)$, a non-commutative probability space with values in $A$, $A\in\mathbf{cAlg}_k$, define the affine subspace $\mathcal{A}_1:=\{a\in\mathcal{A}~|~\phi(a)=1\}$. Then we have the following sequence of pullback cones:
$$
\begin{xy}
  \xymatrix{
  \mathcal{A}_1
  \ar[d]^{R} \ar[rr]^{\mathcal{M}}
  &     &  \Sigma_1(A)\ar[d]^{S}\ar[rr]^{\operatorname{LOG}} &&\Sigma_0(A)\ar[d]^{\mathcal{R}}
  \\
\Sigma_1(A)\ar[rr]^{\mathcal{F}}           &    & \Lambda(A)\ar[rr]^{z\frac{d}{dz}\ln}  & &\mathbb{A}^{\N^*}(A)
}
\end{xy}
$$
with $\mathcal{F}$ as in~(\ref{RScumul}).
\subsection{Witt vectors, $\lambda$-rings and necklace algebras}

The algebraic structures we shall encounter are of the same nature as in~(\ref{V_P_n}) and (\ref{V_Q_n}), despite the fact that they arise in quite different contexts. References for this section are~\cite{DS1989,Hz,I1979,Len,MR}.

Let $\mathbf{cRing}_k$ be the category of commutative unital rings over $k$. Let   
$\W$ be the functor $\W:\mathbf{cAlg}_k\rightarrow\mathbf{cRing}_k$,  with underlying set  $\W(A):=A^{\N^*}$.

\begin{prop}[e.g. \cite{Hz} p. 12]
$(\W(A),+_{\W},\cdot_{\W},(0,0,0,\dots),(1,0,0,0,\dots))$ is a commutative unital ring, called the {\bf ring of Witt vectors}, with multiplicative unit $(1,0,\dots,0,\dots)$. Addition and multiplication are given in terms of universal polynomials $S_{\W}=({S_{\W}}_n)_{n\in\N^*}$, $P_{\W}=({P_{\W}}_n)_{n\in\N^*}$, i.e.
\begin{eqnarray*}
\label{Ws}
({x}+_{\W}{y})_n & := & {S_{\W}}_n({x},{y}), \\
\label{Wm}
({x}\cdot_{\W}{y})_n & := & {P_{\W}}_n({x},{y}),
\end{eqnarray*}
where
${S_{\W}}_n,{P_{\W}}_n\in\Z[x_1,\dots,x_n,y_1,\dots,y_n]$,
and which satisfy:
\begin{itemize}
\item $S_{\W_1}(x_1,y_2)=x_1+y_1$
\item $S_{\W_n}(x_1,\dots,x_n,y_1,\dots,y_n)=S_{\W_n}(y_1,\dots,y_n,x_1,\dots,x_n)$ 
\item $S_{\W_n}(x,y)=x_n+y_n+\tilde{S}_{\W_n}(x_1,\dots,x_{n-1},y_1,\dots,y_{n-1})$ 
\item $P_{\W_1}(x_1,y_1)=x_1\cdot y_1$
\item $P_{\W_n}(x_1,\dots,x_n,y_1,\dots,y_n)=P_{\W_n}(y_1,\dots,y_n,x_1,\dots,x_n)$ 
\end{itemize}
\end{prop}
The polynomials $S_{\W}$ and $P_{\W}$ can be determined by the requirement that the {\bf ghost map} 
\begin{equation*}
w_A:\mathbb{W}(A)\rightarrow A^{\N^*}
\end{equation*}
is a natural isomorphism of unital commutative rings, where $A^{\N^*}$ is given the ring structure of component-wise addition and multiplication. 

The map $w=(w_n)_{n\in\N^*}$,  $w_n(x_1,\dots,x_n)\in\Z[x_1,\dots,x_n]$, composed of the {\bf Witt polynomials}, 
has components 
\begin{eqnarray*}
\label{ghost_map}
w_n&:&\W(A)\rightarrow A,\nonumber\\
{x}&\mapsto& w_n(x_1,\dots,x_n):=\sum_{d|n}d x_d^{n/d},
\end{eqnarray*} 
such that all of them are ring homomorphisms.

The set $\Lambda(A)$ is a commutative ring, the {\bf Grothendieck $\lambda$-ring}, with the addition $+_{\Lambda}$ defined by the usual multiplication of power series, i.e. $+_{\Lambda}:=\cdot$ and with ``zero" $0_{\Lambda}=1$.

The {``multiplication"} $\cdot_{\Lambda}$, is defined as follows: for $f,g\in\Lambda(A)$ consider the formal factorisations 
$f(z)=\prod_{i=1}^{\infty}(1-x_iz)^{-1}$ and $g(z)=\prod_{j=1}^{\infty}(1-y_jz)^{-1}$,
and define 
$$
(1-az)^{-1}\cdot_{\Lambda}(1-bz)^{-1}:=(1-abz)^{-1},
$$ 
which is then bilinearily extended , i.e. one uses
$$
f\cdot_{\Lambda}(g+_{\Lambda}h)=f\cdot_{\Lambda}g+_{\Lambda}f\cdot_{\Lambda}h.
$$
Then for every $A\in\mathbf{cAlg}_k$, the map
\begin{equation*}
z\frac{d}{dz}\ln:(\Lambda(A),+_{\Lambda},\cdot_{\Lambda})\rightarrow (A^{\N^*},\text{$+,\cdot$, pointwise})
\end{equation*}
is an isomorphism of commutative unital rings.

The natural isomorphism $E:\W\rightarrow\Lambda$, the {\bf Artin-Hasse exponential map}, is,  for $A\in\mathbf{cAlg_k}$, given by
$$
E_A:\W(A)\rightarrow\Lambda(A),\quad (x_n)_{n\in\N^*}\mapsto\prod_{n=1}^{\infty}\frac{1}{1-x_nz^n}.
$$

N. Metropolis and G.-C. Rota~\cite{MR}, introduced the 
functor $\operatorname{Nr}:\mathbf{cRing}_k\rightarrow\mathbf{cRing}_k$, called the {\bf necklace algebra}, which has as underlying set $\operatorname{Nr}(A):=A^{\N^*}$.

For $a,b\in\operatorname{Nr}(A)$, addition is given by component-wise addition, and the multiplication $\ast_{\operatorname{MR}}$ is defined as 
\begin{equation*}
\label{}
(a\ast_{\operatorname{MR}}b)_n:=\sum_{\operatorname{lcm}(i,j)=n}\gcd(i,j)\,a_i b_j,
\end{equation*}
where $\operatorname{lcm}(i,j)$ and $\gcd(i,j)$ denote the `least common multiple' and the `greatest common divisor' of $i$ and $j$, respectively. 
The ghost map $g_A:\operatorname{Nr}(A)\rightarrow A^{\N^*}$ is given by
$$
(a_n)_{n\in\N^*}\mapsto\sum_{d|n}d\cdot a_d.
$$

A. Dress and B. Siebeneicher~\cite{DS1989} gave a combinatorial interpretation of the integer Witt ring $\mathbb{W}({\Z})$ by viewing the necklace algebra $\operatorname{Nr}(\Z)$ as the completed {\bf Burnside-Grothendieck ring} $\hat{\Omega}(\mathbf{C})$ of almost finite cyclic sets, where $\mathbf{C}$ is the infinite cyclic group.

\begin{thm}
\label{WittLambda}
For every $A\in\mathbf{cAlg}_k$, the abelian groups $(\W(A),+_{\W})$, $(\operatorname{Nr}(A),+_{\operatorname{Nr}})$, $(\Lambda(A),+_{\Lambda})$, $(\Sigma_0(A),\boxplus)$ and $(\Sigma_1(A),\boxtimes)$ are naturally isomorphic to $(A^{\N^*},+)$, as shown in the diagram below:
\[
\begin{xy}
  \xymatrix{
  \left(\W,+_W\right)\ar[r]^{E}&(\Lambda,\cdot)\ar[dd]^{z\frac{d}{dz}\ln}&\left(\Sigma_1,\boxtimes\right)\ar[l]_S\\
  &&\\
\left(\operatorname{Nr},+\right)\ar[r]^g&(\mathbb{A}^{\N^*},+)&\left(\Sigma_0,\boxplus\right)\ar[l]_{\quad\mathcal{R}} 
     }
\end{xy}
\]
The natural isomorphism $\varphi:(\W,+_{\W})\rightarrow(\Sigma_1,\boxtimes)$,  $\varphi:=E\circ S^{-1}$, is, for $A\in\mathbf{cAlg}_k$, given by
$$
\varphi_A=\operatorname{inv}\left(\frac{z}{1+z}\prod_{n=1}^{\infty}\frac{1}{1-a_nz^n}\right),
$$
where $\operatorname{inv}$ denotes the compositional inverse of power series. Therefore, the equality
$$
\varphi_A(a+_{\W}b)=\varphi_A(a)\boxtimes\varphi_A(b).
$$
holds.

\begin{rem}
In~\cite{Len} C.~Lenart summarises various previously known classical isomorphisms between the different rings. 
\end{rem}
\end{thm}
\subsection{The induced ring structure}
The natural ring structure on $(\mathbb{A}^{\N^*},+,\cdot)$ induces via the $R$- and $S$-transform additional operations.  
For $a,b\in A^{\N^*}$, define
\begin{equation*}
a\boxdot b:=R^{-1}(R(a)\cdot_{\operatorname{H}} R(b)),
\end{equation*}
where $\cdot_{\operatorname{H}}$ is the point-wise Hadamard multiplication,
and for $a,b\in \Sigma_1(A)$, define
\begin{equation*}
a\square\hspace{-0.85 em}\ast b:=S^{-1}(S(a)\cdot_{\Lambda} S(b)).
\end{equation*}

\begin{prop} 
For every $A\in \mathbf{cAlg}_k$, 
$(\Sigma(A),\boxplus,\boxdot)$ and 
$(\Sigma_1(A),\boxtimes, {\square\hspace{-0.6 em}\ast})$ are commutative unital rings  with neutral element $1$ and $S^{-1}(1-z)$, respectively. 

The natural ring isomorphism $\operatorname{LOG}:(\Sigma_1,\boxtimes, {\square\hspace{-0.6 em}\ast})\rightarrow(\Sigma_0,\boxplus,\boxdot)$, is given by 
$$
\operatorname{LOG}:=\mathcal{R}^{-1}\circ(z\frac{d}{dz}\ln)\circ S.
$$
\end{prop}

For compactly supported {\bf freely infinitely-divisible} probability measures on the real line, the structure resulting from $\boxdot$ is a semiring, with a similar statement for ${\square\hspace{-0.65 em}\ast}$\,, as shown in~\cite{FMcK2013b,F2019}. 

\subsection{Faber polynomials and Adams operations}
\label{sec:Faber}
In this section we show the relations Faber polynomials and the Grunsky matrix have with (second order) free probability, Voiculescu's free entropy~\cite{V1998} and integrable systems~\cite{Ta2001,Teo2003}. 

The {\bf Faber polynomials} appear 1903 in complex analysis in approximation theory, see~\cite{McKS,P}. In slightly modified versions they appear later in other contexts. The reason for putting the emphasis on the complex variables approach is that the original Faber polynomials play a profound role in diverse mathematical fields, and in particular in the theory of univalent functions,~see~\cite{McKS,P,Teo2003} and references therein. The explicit connection with Witt vectors reveals new perspectives.

Let $h(z)=1+b_1z+b_2z^2+\cdots\in\Lambda(A)$ with coefficients $b_i\in\C$ and let $g(z):=z h(1/z)$. Then
$$
g(z)=z+\sum_{n=0}^{\infty} b_{n+1}z^{-n}~,
$$
which in complex analytic terms represents the germ of a {\bf univalent function} around the point at infinity, with $g(\infty)=\infty$ and $g'(\infty)=1$. The transformation $g(z)\mapsto\frac{1}{g(1/z)}$ gives a bijection with the locally univalent functions around the origin of the form $z+a_2z^2+a_3z^3+\dots$. 

The {\bf $n$th Faber polynomial} $F_n(w)$ of $g(z)$ is defined by the expansion at infinity
\begin{equation*}
\label{Faber_n(w)}
\ln\frac{g(z)-w}{z}=-\sum_{n=1}^{\infty}\frac{F_n(w)}{n}z^{-n}~.
\end{equation*}
$F_n(w)$ depends on the coefficients $b_1,\dots,b_n$ and the $F_n(w)$ satisfy the recursion relations
$$
F_{n+1}(w)=(w-b_1)F_n(w)-\sum_{k=1}^{n-1}b_{n-k+1}F_k(w)-(n+1)b_{n+1}~.
$$
The {\bf Grunsky matrix} is defined as 
\begin{equation}  
\label{Grunsky}
\ln\frac{g(z)-g(w)}{z-w}=-\sum_{m=1}^{\infty}\sum_{n=1}^{\infty} \beta_{mn} z^{-m}w^{-n}~,
\end{equation}
with the {\bf Grunsky coefficients} $\beta_{mn}$ being symmetric, i.e. $\beta_{mn}=\beta_{nm}$.
The relation between the Faber polynomials and the Grunsky coefficients is given by
\begin{equation*}
\label{FabGrun}
F_n(g(z))=z^n+n\sum_{m=1}^{\infty} \beta_{mn} z^{-m}~.
\end{equation*}
Considering $F_n(0)=:F_n(b_1,\dots,b_n)$ leads to an alternative definition. 
The {\bf $n$th Faber polynomial} $F_n(b_1,\dots, b_n)$ for $n\in\N$ is  defined by 
\begin{equation*}
\label{Faber1}
z\frac{d}{dz}\ln g(z)=\sum_{n=0}^{\infty} F_n(b_1,\dots, b_n) z^{-n}~,
\end{equation*}
with $F_0=1$, $F_1=-b_1$, $F_2=b_1^2-2b_2$, $F_3=-b_1^3+3b_1b_2-3b_3,\dots$. 
  
Equivalently we have, cf.~\cite{Bou},
\begin{equation*}
\label{Bou}
1+b_1z+b_2z^2+b_3z^3+\dots=\exp\left(-\sum_{n=1}^{\infty}\frac{F_n(b_1, b_2,\dots,b_n)}{n}z^n\right)~.
\end{equation*}
\begin{prop}
For a distribution $\mu\in\Sigma_1$, the linearisation of its $S$-transform $S_{\mu}(z)$, is given in terms of the {\bf Faber polynomials} $\{F_n\}_{n\in\N}$ as: 
\begin{equation}
\label{Fab_S}
-z\ln(S_{\mu}(z))=\sum_{n=1}^{\infty}\frac{1}{n}F_n(\mu)z^n. 
\end{equation}
\end{prop}

Alternatively, cf.~\cite{McKS}, the Faber polynomials $F_n(w)$ are given by:
$$
F_n(w)=\det(w1-A_n)~,
$$
where
$$
A_n:=\left(\begin{array}{ccccc}b_1 & 1 & 0 & \cdots & 0 \\ 2b_2 & b_1 & 1 & \ddots & \vdots \\3b_3 & b_2 & \ddots & \ddots & 0 \\\vdots & \vdots & \ddots & b_1 & 1 \\nb_{n} & b_{n-1} & \cdots & b_2 & b_1\end{array}\right)
$$

The {\bf $n$th~Adams operator $\Psi^n$} is given in terms of the {\bf operations $\lambda^j$} by
$$
\Psi^n=\det\left(\begin{array}{ccccc}\lambda^1 & 1 & 0 & \cdots & 0 \\2\lambda^2 & \lambda^1 & 1 & \ddots & \vdots \\3\lambda^3 & \lambda^2 & \ddots & \ddots & 0 \\\vdots & \vdots & \ddots & \lambda^1 & 1 \\n\lambda^n & \lambda^{n-1} & \cdots & \lambda^2 & \lambda^1\end{array}\right)
$$
which corresponds to the determinant formula linking the {\bf power sums} and the {\bf elementary symmetric functions}~\cite{Hz}.

We have the following relation of Faber polynomials with Adams operations.
\begin{prop}
The ghost components of  $\Lambda(k)$ are given by the Faber polynomials $F_n(0)$, $n\in\N$. The $n$th~Adams operator is, up to a sign, equal to the $n$th~Faber polynomial for $w=0$,
$$
\Psi^n=(-1)^nF_n(\lambda^1,\dots,\lambda^n)~,
$$ 
with the coefficient $b_j$ replaced by the operation $\lambda^j$.
\end{prop}

In Voiculescu's approach to the microstate-free entropy~\cite{V1998}, in the one-variable case, he introduced the {\bf difference quotient derivation} $\partial$, as 
\begin{equation*}
g\mapsto\partial g:=\frac{g(x)-g(y)}{x-y}.
\end{equation*}

This motivates to define the {\bf generalised Grunsky matrix} as
\begin{equation*}
g\mapsto\log(\partial g)=\log\frac{g(x)-g(y)}{x-y}.
\end{equation*}
We obtain an identity which equates Voiculescu's free entropy~\cite{V1998} with the free Boson~[\cite{Ta2001} Theorem 3.9. and Corollary 3.12.], via second-order freeness~[\cite{CMSS} Equation~(7)], as
\begin{equation}
\frac{\partial^2}{\partial x\partial y}\log(\partial g)=\sum_{m,n=1}^{\infty}\frac{\partial^2\log(\tau)}{\partial t_m\partial t_n}x^{-m-1}y^{-n-1}=G(x,y),
\end{equation}
where $\tau$ is the $\tau$-function of a smooth Jordan contour $C$ as in~\cite{Ta2001}, cf. also~\cite{Teo2003}, and $G(x,y)$ the {\bf second order Cauchy transform} in the case of vanishing {\bf second order free cumulants}, as for Gaussian or Wishart random matrices~\cite{CMSS}.

The {\bf generalised energy-moment tensor} or {\bf generalised Schwarzian derivative}, cf. for conformal field theory~\cite{KNTY} and for the connection with {\bf replicable functions}~\cite{McKS}, is then 
\begin{equation*}
g\mapsto\lim_{z\to w}d_zd_w\log(\partial g)=\lim_{z\to w}d_zd_w\log\frac{g(x)-g(y)}{x-y}.
\end{equation*}

\section{Multiplicative genera and Fock spaces}
\subsection{The full Fock space and canonical random variables}
Let us recall from~\cite{Haa1997,NS,V1986,VDN} some facts which are relevant in connection with multiplicative genera but also with Faber polynomials and the Grunsky matrix.

For $\mathcal{H}$ a $\C$-Hilbert space (finite or infinite), the {\bf full Fock space over $\mathcal{H}$} is 
$$
\mathfrak{F}(\mathcal{H}):=T(\mathcal{H})=\bigoplus_{n=0}^{\infty}\mathcal{H}^{\otimes n}~\cong\C\Omega\oplus\bigoplus_{n=1}^{\infty}\mathcal{H}^{\otimes n}~,
$$
with $\Omega$ the {\bf vacuum vector}, a distinguished unit vector, and $B(\mathfrak{F}(\mathcal{H}))$ the algebra of bounded linear operators. 
$\mathfrak{F}(\mathcal{H})$ is a Hilbert space with the different summands $\mathcal{H}^{\otimes n}$ being pairwise orthogonal, and with
$$
\langle x_1\otimes\cdots\otimes x_n, y_1\otimes\cdots\otimes y_n\rangle=\langle x_1, y_1\rangle\cdots\langle x_n, y_n\rangle~.
$$
For $T\in B\left(\mathfrak{F}(\mathcal{H})\right)$, the {\bf vacuum expectation value} $\tau_{\Omega}$ is given by 
\begin{equation*}
\tau_{\Omega}(T):=\langle T(\Omega),\Omega\rangle~,
\end{equation*}
which induces the {\bf vacuum state} $\tau_{\Omega}:B\left(\mathfrak{F}(\mathcal{H})\right)\rightarrow\C$, such that $(B(\mathfrak{F}(\mathcal{H})),\tau_{\Omega})$ becomes a $C^*$-probability space. 

For $h\in\mathcal{H}$, the associated {\bf left-creation operator} $l(h):\mathfrak{F}(\mathcal{H})\rightarrow \mathfrak{F}(\mathcal{H})$ is defined as
\begin{eqnarray*}
l(h)\Omega & := &h, \\
l(h) x_1\otimes\cdots\otimes x_n & := & h\otimes x_1\otimes\cdots\otimes x_n,
\end{eqnarray*}
with the {\bf adjoint operator} $l^*(h)$ satisfying
\begin{eqnarray*}
l^*(h)\Omega & = &0 \\
l^*(h)x_1\otimes\cdots\otimes x_n&=&\begin{cases}
      & \langle h,x_1\rangle\,\Omega,\qquad\qquad\qquad\qquad\;\,\text{for $n=1$}, \\
      & \langle h, x_1\rangle\, x_2\otimes\cdots\otimes x_n\qquad\qquad \text{for $n\geq 2$}.
\end{cases}
\end{eqnarray*}
For $e_1,\dots,e_n$, pairwise orthonormal vectors in $\mathcal{H}$, let $l_i:=l(e_1)$, for $i=1,\dots, n$.

Voiculescu~[\cite{VDN} Lemma and Definition 3.2.2.] gave a functional analytic representations of the $R$-transform in terms of ``canonical random variables", as did U. Haagerup~[\cite{Haa1997} Theorem 2.2. and Theorem 2.3], who additionally found such a representation for the $S$-transform, which we both recall now. One should also note the parallels with the derivation of the Faber polynomials in Section~\ref{sec:Faber}. 

\begin{thm}[Fock space models for the $R$- and $S$-transform]
\label{Fock-canonical}
For every $\mu\in\Sigma$, with moment series $\mathcal{M}_{\mu}(z)$, moments $(m_n(\mu))_{n\in\N^*}$ and free cumulants $(\kappa_n(\mu))_{n\in\N^*}$, define for the {\bf canonical random variables}
\begin{eqnarray*}
\label{crv}
T_{\mu}&:=&l^*_1+\sum_{n=0}^{\infty}\kappa_{n+1}(\mu){l_1}^n=l^*_1+\mathcal{R}_{\mu}(l_1),\\
T_{\mu}^*&:=&l_1+\sum_{n=0}^{\infty}m_{n+1}(\mu){(l^*_1)}^n=l_1+\mathcal{M}_{\mu}(l_1^*).\\
\end{eqnarray*}
Then the identities 
\begin{equation}
\mathcal{M}_{l_1^*+\mathcal{R}_{\mu}(l_1)}(z)=\mathcal{M}_{\mu}(z)
=\mathcal{R}_{l_1+\mathcal{M}_{\mu}(l_1^*)}(z)
\end{equation}
hold. 

For every $\nu\in\Sigma^{\times}$, with moment series $\mathcal{M}_{\nu}(z)$, define the operator
\begin{equation*}
\label{H_S-trafo}
a_{\nu}:=(\id+l_1)\mathcal{M}_{\nu}(l_1^*).
\end{equation*}
Then the identities  
\begin{equation}
\label{Haa_S-trafo}
S_{a_{\nu}}(z)=\frac{1}{\mathcal{M}_{\nu}(z)}
\end{equation}
holds.
\end{thm}

\subsection{Multiplicative genera}
Previously, in~\cite{FMcK2011} we described relations between free probability and complex cobordism, which we shall now extend further.
As references we use~\cite{A,BMN,H1956,KSU}.

Let $MU^*$ denote the complex cobordism ring and $A\in\mathbf{cAlg}_{\Q}$. 

A {\bf genus} with values in $A$, 
is a ring homomorphism 
$$
\varphi:MU^*\otimes_{\mathbb
Z}\mathbb{Q}\rightarrow A,\quad\text{with}\quad \varphi(1_{MU^*\otimes_{\mathbb Z}\mathbb{Q}})=1_A.
$$

For $\C{\mathbb P}^n$, the $n$-dimensional complex projective space, with $\C{\mathbb P}^0:=1$, according to A. Mishchenko~\cite{BMN}, the {\bf logarithm of the genus $\varphi$}, is defined as the formal power series
\begin{equation}
\label{log_cob}
\log_{\varphi}(z):=\sum_{n=1}^{\infty}\frac{1}{n}\varphi(\C{\mathbb
P}^{n-1})z^n\in\operatorname{Aut}_+(\mathcal{O}_{A}),
\end{equation}
with inverse $\log^{-1}_{\varphi}(z)$.

According to Novikov~\cite{BMN}, every genus $\varphi$ determines uniquely a {\bf Hirzebruch multiplicative sequence} $\{K_n\}$, cf.~\cite{H1956}, via the associated {\bf characteristic power series} $Q_{\varphi}(z)\in\Lambda(A)$, through the identities 
\begin{equation} 
\label{NovHir}
K(1+z)=Q_{\varphi}(z):=\frac{z}{\log^{-1}_{\varphi}(z)}.
\end{equation}
Further, over $\Q$, all Hirzebruch genera yield formal group laws which are strongly isomorphic to the additive group law, cf.~\cite{BMN}. 

Now the following bijections, as shown in the commutative diagrams below, hold. For complex cobordism we have
\begin{equation*}
\label{Genus_square}
\begin{xy}
  \xymatrix{
    \Hom_{\mathbf{cRing}_{\Q}}(MU^*\otimes_{\Z} \Q,A)
  \ar[d]^{\log} \ar[rr]^{\qquad\qquad\text{$Q$-series}}
  &     &  \Lambda(A) 
  \\
\operatorname{Aut}_+(\mathcal{O}_{A})\ar[urr]_{z/\operatorname{log}^{-1}}           &    & }   
\end{xy}
\end{equation*}
and for free probability
\begin{equation*}
\label{R_S_trafo_square}
\begin{xy}
  \xymatrix{
    \Sigma_1(A)=\Hom_{\Q,\bullet,1}(A[x],A)
  \ar[d]^{R} \ar[rr]^{\qquad\qquad S}
  &     &  \Lambda(A)
 \\
\operatorname{Aut}_+(\mathcal{O}_{A})\ar[urr]_{R^{-1}/z}          &    & }
\end{xy}
\end{equation*}
holds.

Haagerup's representation of the $S$-transform~(\ref{Haa_S-trafo}) induces a natural  operator-valued representation of complex genera. 
\begin{prop}
\label{C*^-genera}
There exists a lift $\sigma:\Hom_{\mathbf{cRing}}(MU^*\otimes_{\Z}\Q,\C)\rightarrow (B(\mathfrak{F}(H)),\tau_{\Omega})$, such that the diagram
$$
\begin{xy}
  \xymatrix{  &     &  (B(\mathfrak{F}(\mathcal{H}),\tau_{\Omega})\ar[d]^S  \\
 \Hom_{\mathbf{cRing}}(MU^*\otimes_{\Z} \Q,\C)\ar[urr]^{\sigma}\ar[rr]^{\qquad Q}            &    & \Lambda(k)  }
\end{xy}
$$
commutes. For a genus $\varphi$, define 
\begin{equation*}
a_{\varphi}:=({ 1}+\ell)\frac{1}{Q_{\varphi}}(\ell^*)\in B(\mathfrak{F}(\mathcal{H})),
\end{equation*}
with $Q_{\varphi}$ the characteristic power series associated to $\varphi$, from which we obtain the identity:
\begin{equation*}
S_{a_{\varphi}}(z)=Q_{\varphi}(z).
\end{equation*}
\end{prop}

From the bijection 
$$
\Hom_{\mathbf{cAlg}_{\Q}}(MU^*\otimes_{\Z}\Q,A)\rightarrow\Sigma_1(A),
$$ 
we obtain a relation, already deduced in~\cite{FMcK2011}, between a genus $\varphi$ and a distribution $\tilde{\mu}$, namely
\begin{equation}
\label{relation_R_genus}
R_{\tilde{\mu}}(z)=\log_{\varphi}(z).
\end{equation} 
Then the identities~(\ref{R-cumul}), (\ref{RScumul}) and~(\ref{log_cob}) yield
\begin{equation*}
\kappa_n(\tilde{\mu})=\frac{1}{n}\varphi[\C\P^{n-1}].
\end{equation*}

\subsection{Genera and probability measures}

In order to relate genera with probability measures, a gauge choice is needed. This is permitted, as cumulants are defined only up to a linear transformation, given by an infinite lower-triangular matrix.

Namely, the linearising map from ``moments" $(x_n)_{n\in\N^*}$ to ``cumulants" $(y_n)_{n\in\N^*}$  is given by universal polynomials $(p_n)_{n\in\N^*}$ such that $y_n=p_n(x_1,\dots,x_n)$. Therefore every linear automorphism of the additive group $(\mathbb{A}^{\N^*},+)$ of cumulants which preservs the causal structure, is given by a projective limit of invertible lower-triangular matrices.

So, let us reformulate equation~(\ref{relation_R_genus}) in terms of the invariant differential~(\ref{inv_diff}) as follows. 
Consider the map $\varphi\mapsto\mu$, given by
\begin{equation*}
\mathcal{R}_{\mu}(z)dz=d\log_{\varphi}(z),
\end{equation*}
from which the identity 
\begin{equation*}
\kappa_n(\mu)=\varphi[\C\P^{n-1}]
\end{equation*}
follows.

Now, T.~Honda's~\cite{Ho1968} definition of a zeta-function for group laws, motivates the following formal $L$-series for measures.
\begin{df}
For $\mu\in\Sigma_1$, with free cumulants $(\kappa_n(\mu))_{n\in\N^*}$,  
the {\bf $\zeta$-function} of the distribution $\mu$ is given by:
\begin{equation}
\label{measure_zeta}
\zeta_{\mu}(s):=\sum_{n=1}^{\infty}\kappa_{n+1}(\mu)n^{-s}.
\end{equation}
\end{df}

Let us discuss now three particular examples. First, we recall from~\cite{NS,VDN} the following facts.
\begin{df}
Let $\lambda\in\R_+$ and $\alpha\in\R$. The {\bf free Poisson} distribution with {\bf rate}, $\lambda$, and {\bf jump size}, $\alpha$, is given by the limit in distribution
\begin{equation*}
\nu_{\infty,\lambda,\alpha}:=\lim_{N\to\infty}\left(\left(1-\frac{\lambda}{N}\right)\delta_0+\frac{\lambda}{N}\delta_{\alpha}\right)^{\boxplus N}.
\end{equation*}
\end{df}
\begin{df}
The {\bf Marchenko-Pastur} measure, with parameters $\lambda\in\R^*_+$ and $\alpha\in\R$, is the probability measure $\mu_{\operatorname{MP},\lambda,\alpha}$ on $\R$, with support the interval
\begin{equation*}
[\alpha(1-\sqrt{\lambda})^2,\alpha(1+\sqrt{\lambda})^2],
\end{equation*}
and which is given by 
\begin{equation*}
\mu_{\operatorname{MP},\lambda,\alpha}:=\begin{cases}
      & (1-\lambda)\delta_0+\lambda\nu\quad \text{for $0\leq\lambda\leq1$ }, \\
      & \nu\quad\text{for $\lambda>1$},
\end{cases}
\end{equation*}
where $\nu$ is has the density
\begin{equation*}
d\nu(x)=\frac{1}{2\pi\alpha x}\sqrt{4\lambda\alpha^2-\left(x-\alpha(1+\lambda)\right)^2}dx.
\end{equation*}
\end{df}
The following statement summarises the content of~[\cite{VDN} (c) pp. 34-35]  and~[\cite{NS} pp. 218, 380].
\begin{prop}
For $\lambda\in\R^*_+$ and $\alpha\in\R$, the equality 
\begin{equation*}
\nu_{\infty,\lambda,\alpha}=\mu_{\operatorname{MP},\lambda,\alpha}
\end{equation*}
holds in distribution, and the corresponding 
free cumulants are $\kappa_n(\nu_{\infty,\lambda,\alpha})=\lambda\alpha^n$, $n\in\N^*$. Its $R$-transform is
\begin{equation*}
R_{\nu_{\infty,\lambda,\alpha}}(z)=\lambda\alpha\frac{z}{1-\alpha z}.
\end{equation*}
\end{prop}
\begin{df}
The {\bf semi-circular} distribution with centre $a$ and radius $r$, is the probability measure with support the compact interval $[-r+a,a+r]$ and with density 
$$
f(x)=
\begin{cases}
  \frac{2}{\pi r^2}\sqrt{r^2-x^2}    & \text{for $-r+a<x<a+r$}, \\
      & \text{$0$ otherwise}.
\end{cases}
$$ 
\end{df}

\begin{prop}
\label{semi_ciruclar_genus}
Let $\gamma_{1,2}$ be the semi-circular distribution centred at $1$ and of radius $2$, and $R_{\gamma_{1,2}}(z)=z+z^2$ its $R$-transform. Then the $S$-transform of $\gamma_{1,2}$ is equal to  
\begin{equation*}
S_{\gamma_{1,2}}(z)=\frac{1}{2z}\left(\sqrt{4z+1}-1\right)=\sum_{n=1}^{\infty} (-1^{n-1})C_n z^{n-1},
\end{equation*}
where $C_n$ is the $n$th {\bf Catalan number}.
\end{prop}

Let us discuss the following examples and summarise the corresponding results in a table.
\begin{enumerate}
\item The {\bf trivial genus}, with characteristic power series $Q(z)=1$, corresponds to the Dirac distribution at $1$.
\item The {\bf Todd genus}, with logarithm $-\log(1-x)$, corresponds to the free Poisson distribution with parameters $\lambda=\alpha=1$, i.e. 
$$
z\frac{d}{dz}\log_{Td}(z)=R_{\nu_{\infty,1,1}}(z).
$$
\item For the {\bf semi-circular} distribution, $\gamma_{1,2}$, the corresponding characteristic power series $Q(z)$ is equal to $z/(-1+\sqrt{2z+1})$.
\end{enumerate}

\[
\begin{tabular}{c|c|c|c} 
Distribution: & Name: & multiplicative genus & name
\\\hline 
$\delta_1$ & Dirac & 1 & trivial genus\\
$\nu_{\infty,1,1}$ & free Poisson & $\frac{z}{1-e^{-z}}$ &  Todd genus\\
$\gamma_{1,2}$ & semi-circular & $\frac{z}{-1+\sqrt{2z+1}}$ & ???
\end{tabular}
\]
In the {\bf Planck's radiation law}, cf.~\cite{Wiki},  according to the {\bf Bose-Einstein statistics}, the {\bf average energy of quantised modes} is given by the formula
\begin{equation}
\label{A_Energy}
E(\nu,T)=\frac{\hbar\nu}{e^{\frac{\hbar\nu}{k_B\cdot T}}-1},
\end{equation}
where $\hbar$ is Planck's Wirkungsquantum, $\nu$ frequency, $T$ temperature and $k_B$ the Boltzmann constant. By performing in~(\ref{A_Energy}) the formal substitution $\nu\mapsto-\nu$ and by setting $T=k_B=\hbar=1$, we obtain
\begin{equation*}
E(-\nu,1)=\varphi_{Td}(\nu).
\end{equation*} 
\section{Physical applications}
We continue the discussion of free probability from a conformal field theoretic perspective, as started in~\cite{FMcK2011}. 
Novel aspects will include S. Kerov's~\cite{K} asymptotic representation theory of the symmetric group. Relations of the Hopf algebra $\operatorname{Symm}$ with the Fock space have previously appeared in~\cite{LTh,Mo}, and as we show, these results are natural from the point of view of Witt vectors and the Sato-Segal-Wilson Grassmannian. Background material covering the topics involved can be found in~\cite{AGMV,KSU,KNTY,MJD,Mo}.
\subsection{Fermion Fock space}
Generally speaking, the {\bf Boson-Fermion correspondence}~\cite{AGMV,KNTY,MJD} is a rich source of algebraic identities, with analogs in non-commutative probability theory. 

\begin{df}[\cite{KNTY,MJD}]
A {\bf Maya diagram} consists of a subset $M\subset\Z+1/2$ such that $M\cap(\Z +1/2)_{>0}$ and $\complement M\cap(\Z+1/2)_{<0}$ are finite sets, where $(\Z+1/2)_{>0}:=\{u\in\Z+1/2~|~u>0\}$ etc., and $\complement M$ is the complement of $M$ in $\Z+1/2$. 

The {\bf charge} or {\bf Euler characteristic} of a Maya diagram $M$, is the integer number
$$
\chi(M):=\#(M\cap(\Z+1/2)_{>0})-\#(\complement M\cap(\Z+1/2)_{<0}).
$$
\end{df}
We denote the set of all Maya diagrams by $\mathcal{M}$
and the subset of diagrams of {\bf charge $p$}, $p\in\Z$, by $\mathcal{M}_p$. The {\bf Fermion Fock} space is the $\C$-vector space spanned by all Maya diagrams.

\begin{figure}[htbp]
\begin{center}
\includegraphics[scale=0.6]{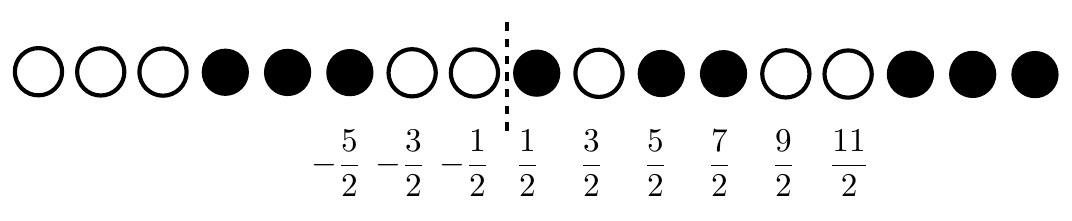}
\caption{Maya diagram of charge $0$ corresponding to the state $|-\frac{9}{2},-\frac{7}{2},-\frac{5}{2};\frac{3}{2}, \frac{9}{2}, \frac{11}{2}\rangle$.}
\label{FMayaDiagram}
\end{center}
\end{figure}

A pictorial definition of a Maya diagram, as given in~\cite{KNTY,MJD}, consists of an infinite sequence of black $\bullet$ and white $\circ$ dots ({\bf Go stones}) positioned at the half integers $\Z+1/2$, such that on the positive axis there are only finitely many white $\circ$ and on the negative axis only finitely many black $\bullet$ dots, as shown in ~Figure~\ref{FMayaDiagram}.

The {\bf Dirac vacuum} $|0\rangle$ corresponds to the sequence with all white stones lined up on the negative axis and the black stones on the positive axis. 

Any Maya diagram, up to a sign, can be obtained by successively applying {\bf fermionic creation operators} $\psi_m,\psi^*_n$, $m,n<0$, to the Dirac vacuum, i.e.
$$
|n_1\dots n_s;m_1\dots m_{r}\rangle:=\psi_{m_1}\cdots\psi_{m_r}\psi^*_{n_1}\cdots\psi^*_{n_s}|0\rangle,
$$
for $m_1<\cdots<m_r<0$ and $n_1<\cdots<n_s<0$. The first $r$ entries correspond to {\bf particle states} and the last $s$ entries to {\bf anti-particle states}.
The bijection between Maya diagrams and charged Young diagrams, cf. e.g. ~[\cite{KNTY} (1.8) Lemma] permits us to give fermionic interpretations of several (free) probabilistic results.

\subsection{Asymptotic representation theory}
S. Kerov~\cite{K} defined maps from the set $\Y_0$ of {\bf Young diagrams} of charge zero to the space of probability measures on the real line, i.e.
\begin{equation}
\label{Kerov}
\operatorname{Kerov}:\Y_0\rightarrow\text{$\P$-measures on $\R$},
\end{equation}
by associating to $\lambda\in\Y_0$, the {\bf transition} $\hat{\omega}_{\lambda}$ and the {\bf co-transition} $\check{\omega}_{\lambda}$ measure, respectively. 

These are 
defined as follows, cf.~\cite{HO,K,KLM}:
\begin{eqnarray*}
\hat{\omega}_{\lambda}&:=&\sum^r_{i=1}\hat{q}_{\lambda}(i)\,\delta_{x_i},\\
\check{\omega}_{\lambda}&:=&\sum_{i=1}^{r-1}\check{q}_{\lambda}(i)\,\delta_{y_j},
\end{eqnarray*}
where $x_1<y_1\cdots<x_{r-1}<y_{r-1}<x_r$ is an {\bf interlacing sequence} of integers $x_i,y_j$, and $\delta_{x_i},\delta_{y_j}$ Dirac measures supported at $x_i$ and $y_j$, respectively. 

The weights
\begin{eqnarray*}
\hat{q}_{\lambda}(i)&:=&\prod_{j=1}^{i-1}\frac{x_i-y_j}{x_i-x_j}\prod^r_{l=i+1}\frac{x_i-y_{l-1}}{x_i-x_l},\\
\check{q}_{\lambda}(i)&:=&\frac{(x_r-y_i)(y_i-x_1)}{\sum_{j<l}(y_j-x_j)(x_l-y_{l-1})}\prod_{j=1}^{i-1}\frac{y_i-y_{j+1}}{y_i-y_j}\prod^{r}_{l=i+1}\frac{y_i-x_{l}}{y_i-y_l},
\end{eqnarray*}
correspond to the possibilities to add or remove from the $i$th row a single cell; for details cf.~[\cite{HO} pp. 253] or~\cite{B1998,K,KLM}.

Kerov~\cite{K} and Ph.~Biane~\cite{B1998} proved that the limit shape $\Omega$ of rescaled Young diagrams with respect to the {\bf Plancherel measure} is equal to, cf.~[\cite{HO} pp. 267],
$$
\Omega(x)=\begin{cases}
     \frac{2}{\pi}(x\arcsin(\frac{x}{2})+\sqrt{4-x^2}), & |x|\leq2, \\
      |x|,& |x|>2.
\end{cases}
$$
The associated transition probability measure $\hat{\omega}_{\Omega}$ is the {\bf Wigner semi-circular} law
$$
\frac{d\hat{{\omega}}_{\Omega}}{dx}=\frac{1}{2\pi}\sqrt{4-x^2},
$$
which corresponds exactly to the distribution from Proposition~\ref{semi_ciruclar_genus}, when translated by $1$.

\subsection{Measures, Grassmannian and $\tau$-function}
Katsura, Shimizu and Ueno~\cite{KSU} established a link between CFT and the complex cobordism ring which we extend now further to free probability theory.

Consider the complex Hilbert space of Laurent series $H:=\{\sum_{n\in\Z} a_n z^n~| ~a_n\in\C\}$ with the polarisation $H=H_+\oplus H_-$, and with bases $\{z^{n}| n\geq0\}$ for  $H_+$ and  $\{z^{m}| m<0\}$ for $H_{-}$, respectively.

The associated {\bf (Sato)-Segal-Wilson Grassmannian}~\cite{SW1985}, $\operatorname{Gr}(H)$, is the set of subspaces $W\subset H$, such that the orthogonal projection $\pr_+:W\rightarrow H_+$ is a Fredholm operator and $\pr_{-}:W\rightarrow H_{-}$ a Hilbert-Schmidt operator. The {\bf big cell} $\operatorname{Gr}_0(H)$ is $\{W\in\operatorname{Gr}(H)~|~\text{$\pr_+:W\rightarrow H_+$ isomorphism} \}$.

A. Kirillov and D.~Yuriev~\cite{KY1988} showed that every Grunsky matrix~(\ref{Grunsky}) defines a subspace $W\in\operatorname{Gr}(H)$, spanned by elements of the form
$$
w_n:=z^{-n}+\sum^{\infty}_{m=1} B_{nm}z^{m-1}.
$$
Physically, as we observe, this corresponds to a {\bf Bogoliubov transformation} of the Dirac vacuum~\cite{AGMV}, which is represented in terms of fermionic operators  as:
$$
|B(W)\rangle=\exp\left(-\sum_{n,m=1}^{\infty}B_{nm}\psi_{-m+1/2} \psi^*_{-n+1/2}\right)|0\rangle.
$$

Every complex {\bf multiplicative sequence} $K = \{ K_{n} (p_{1}, \dots , p_{n}) \}_{n \in \mathbb{N}}$, as introduced by Hirzebruch~\cite{H1956}, is a sequence, for $\deg(p_i)=j$, of homogeneous polynomials with complex coefficients which gives rise to a monoid isomorphism 
$$
K: \Lambda ( \mathbb{C} ) \ni
1 + \sum_{j=1}^{\infty} p_{j} z^{j}
\mapsto
\sum_{j=0}^{\infty} K_{j} ( p_{1}, \dots , p_{j} ) z^{j}
\in \mathbb{C} [p_1,p_2,p_3,\dots][\![ z ]\!],
$$
with
$K_{0} := 1$.

Let $b_{n} \in \mathbb{C}$ be 
the coefficient of $p_{n}$
in the polynomial 
$K_{n} ( p_{1}, \dots , p_{n} )$. 
Then according to [\cite{KSU} Theorem~5.3, 1], the expression 
$$
\tau_{K} ( \mathbf{t} )
:=
\exp \Big(
	\sum_{n=1}^{\infty}
	(-1)^{n+1} b_{n} t_{n}
\Big)\in\hat{\mathcal{H}}_0(\C),
$$
with $\mathbf{t} = ( t_{1}, t_{2},t_{3}, \dots )$ and $\hat{\mathcal{H}}_0(\C)$, the completion of the {\bf Boson Fock space} $\mathcal{H}_0$ of {\bf central charge} zero, is a {\bf tau-function of the KP hierarchy}, cf.~\cite{KNTY,MJD}. 
In particular, the coefficients $b_j$ can be expressed as polynomials in the free cumulants $\kappa_i$, i.e. 
$$
b_j=b_j(\kappa_2,\dots,\kappa_j).
$$

By [\cite{KSU} Theorem 5.3] there exists an injective map from the ring of Witt vectors (Lambda ring) into the big cell of the Grassmannian:
\begin{eqnarray}
\label{KSU}
\operatorname{KSU}:\Lambda&\rightarrow&\operatorname{Gr}_0,\\\nonumber
a&\mapsto&W_a,
\end{eqnarray}
where $a_K=(a_1(K),a_2(K),a_3(K),\dots)$ is the characteristic power series corresponding to a multiplicative sequence $K$, as in equation~(\ref{NovHir}).
In summary, we obtain 
\begin{thm}
Every distribution with mean zero (probability measure with mean zero on the real line) can be embedded into the big cell of the Sato-Segal-Wilson Grassmannian, as shown in the diagram:
\[
\begin{xy}
  \xymatrix{
  \Y_0\ar[r]^{\operatorname{Kerov}}&\Sigma_0\ar[r]^{\operatorname{EXP}}&\Sigma_1\ar[r]^S&\Lambda\ar[r]^{\operatorname{KSU}}&\operatorname{Gr}_0
  } 
\end{xy}
\]
with the injection $\operatorname{Kerov}$~(\ref{Kerov}), bijection $\operatorname{EXP}$, which is the inverse to~(\ref{LOG}), bijection $S$~(\ref{S_trafo}) and injection $\operatorname{KSU}$~(\ref{KSU}).
\end{thm}

The Virasoro algebra of central charge one, as part of a bigger algebra,  acts on the determinant line bundle $\Det$ over $\operatorname{Gr}_0(H)$, cf.~\cite{KNTY,KY1988,MJD}, and hence by our embedding on the set of probability measures. 

First, this relates to the work of, H.~Kvinge, A.~Licata and S.~Mitchell~\cite{KLM} and A.~Lascoux and JY.~Thibon~\cite{LTh} and others naturally, who studied the possible links {\bf Khovanov's Heisenberg category} has with {\bf Jones' planar algebras} and free probability or the action of vertex operators on (shifted) symmetric functions, and in particular polynomial representations of $\operatorname{GL}_n$. From the infinite dimensional perspective we have considered, either Segal and Wilson's~\cite{SW1985} functional analytic restricted general linear group $\operatorname{GL}_{1}(H)^0$ and its $\C^*$-central extension $\hat{\operatorname{GL}}_{1}$ or the Kyoto School's~[\cite{MJD} pp. 53] vertex operator representation of the infinite Lie algebra $\mathfrak{gl}(\infty)$ and its corresponding Lie group $\mathbf{G}$ have to appear.

Second, the {\bf Krichever embedding}~\cite{AGMV,KNTY,KY1988,SW1985}, $(C,L,p,z,\varphi)\mapsto W\in\operatorname{Gr}(H)$, of the datum consisting of a compact Riemann surface $C$ (a curve) of genus $g$, a complex line bundle $L$ over $C$, $p\in C$ a marked point, $z$ a local parameter around $p$, and $\varphi:L|_{\overline{U}_p}\rightarrow\overline{\D}\times\C$ a local trivialisation over the closure of a neighbourhood $U_p$ of $p$, where $\D$ is the unit disc, defines a point $W$ in the Sato-Segal-Wilson Grassmannian. If the {\bf Euler characteristic}, $\chi(L):=\dim_{\C} H^0(C;L)-\dim_{\C} H^1(C;L)$, is equal to $1$, then the corresponding subspace $W$ belongs to the big cell $\operatorname{Gr}_0(H)$~[\cite{SW1985} Proposition 6.1.]. Therefore we obtain for the moduli space $\hat{\mathcal{M}}_{g,1}$ of triples $(C,p,z)$, with $(C,p)\in\mathcal{M}_{g,1}$, for $\mathcal{M}_{g,1}$ the moduli space of Riemann surfaces $C$ of genus $g$ with $p\in C$, a marked point, and $z$ a formal coordinate at $p$, the fibred product, $\hat{\mathcal{M}}_g\times_{\operatorname{Gr}_0}\Lambda$, cf.~[\cite{KNTY} (2.28) Theorem.]. So, it would be interesting to characterise those curves which arise from complex genera / probability measures and to describe the relation with the genus expansion of random matrices. 

Third, as outlined in the work with T.~Amaba~\cite{AF2019b}, besides the known relations between random matrices, integrable hierarchies and growth models, cf. e.g.~\cite{Ta2001,Teo2003}, there exists, as we show, a link with higher order free probability~\cite{CMSS}. The dictionary we started to compile, gives yet another direction, in particular in connection with B.~Eynard's ``Topological Recursion"~\cite{E2016}, to be further studied. 

\subsection*{Acknowledgements}
Both authors thank the MPIfM in Bonn for the opportunity to collaborate on this project on a previous occasion, and for its hospitality. We both thank the referee for the helpful and constructive comments. R.F. was partially supported by the ERC advanced grant ``Noncommutative distributions in free probability".

Authors addresses

Roland Friedrich,\\ 
Fachrichtung Mathematik,\\
Universität des Saarlandes,\\
66123 Saarbrücken,\\ 
Germany\\ 
friedrich@math.uni-sb.de

John McKay,\\ 
Dept. of Mathematics,\\ 
Concordia University,\\ 
Montreal,\\
Canada\\
jmckay1939@icloud.com


\begin{thebibliography}{00}

\bibitem{A} 
J.F. Adams, {\it Stable Homotopy and Generalised Homology}, Chicago Lectures in  Mathematics, The University of Chicago Press, Chicago and London, (1974)

\bibitem{AGMV} I. Alvarez-Gaumé, C. Gomez, G. More, C. Vafa, {\it Strings in the Operator Formalism}, Nucl. Phys. {\bf B} 303, North-Holland, (1988)

\bibitem{AF2019b}
T. Amaba, R. Friedrich, {\it Modulus of continuity of controlled Loewner-Kufarev equations and random matrices}, arXiv (2019)

\bibitem{AA2016} M. Anshelevich , O. Arizmendi {\it The Exponential Map in Non-commutative Probability}, International Mathematics Research Notices, Vol. 2017, No. {\bf 17}, pp. 5302–5342, (2017)

\bibitem{B1998} P.Biane, {\it Representations of symmetric groups and free probability}, Adv. Math. {\bf 138} p. 126–181, (1998).

\bibitem{Bou} A. Bouali, {\it Faber polynomials, Cayley–Hamilton equation and Newton symmetric functions}, Bull. Sci. math. {\bf 130}, 49-70 (2006)

\bibitem{BMN}
V.M. Bukhshtaber, A.S. Mishchenko, S.P. Novikov {\it Formal Groups and their Role in the Apparatus of Algebraic Topology}, (1970)

\bibitem{C2016} G. Cébron, {\it Matricial Model for the Free Multiplicative Convolution}, The Annals of Probability, Vol. {\bf 44}, No. 4, 2427–2478 (2016)

\bibitem{CMSS} B. Collins, J. Mingo, P. $\rm\acute{S}$niady, R. Speicher, {\it Second Order Freeness and Fluctuations of Random Matrices III. Higher Order Freeness and Free Cumulants}, Documenta Math. 12, 1–70, (2007) 

\bibitem{D1972} M. Demazure, {\it Lectures on $p$-divisible groups}, Lecture Notes in Mathematics, Vol. {\bf 302}, Springer, (1972).

\bibitem{DS1989} A. Dress, Ch. Siebeneicher, 
{\it The Burnside Ring of the Infinite Cyclic Group and Its Relations to the Necklace Algebra, $\lambda$-Rings, and the Universal Ring of Witt Vectors}, Adv. in Mathematics {\bf 78}, 1-41 (1989)

\bibitem{DPT} G. Drummond-Cole, J.-S. Park, J. Terilla, {\it Homotopy probability theory I.}, J. Homotopy Relat. Struct. 10, no. {\bf 3}, (2015)

\bibitem{KP2015}  K. Ebrahimi-Fard, F. Patras {\it Cumulants, free cumulants and half-shuffles},
Proc. R. Soc. A:  471, Issue: 2176. (2015)

\bibitem{E2016} B. Eynard, {Counting Surfaces}, CRM Aisenstadt Chair Lectures, Progress in Mathematical Physics, CRM Birkhäuser, (2016)

\bibitem{FrdMcK}
D. Ford, J. McKay, {\it Monstrous Moonshine--two footnotes}, Third spring conference 
``Modular Forms and Related Topics", Hamamatsu, 56-61, (2004)

\bibitem{FM} A. Frabetti, D. Manchon, {\it Five interpretations of Faà di Bruno’s formula, in ``Faà di Bruno Hopf algebras, Dyson-Schwinger equations, and Lie-Butcher series"}, pp. 91–147, IRMA Lect. Math. Theor. Phys., 21, European Mathematical Society, (2015).

\bibitem{FB} E. Frenkel, D. Ben-Zvi, {\it Vertex Algebras and Algebraic Curves}, Second Edition, Mathematical Surveys and Monographs, Volume: {\bf 88}, 400 pp, AMS,
(2004)
\bibitem{FMcK2011}
R. Friedrich, J. McKay, {\it Free Probability Theory and Complex Cobordism},
C. R. Math. Rep. Acad. Sci. Canada Vol. {\bf 33} (4), pp. 116-122 (2011)
\bibitem{FMcK2013} R. Friedrich, J. McKay, {\it Almost Commutative Probability Theory}, arXiv (2013)
\bibitem{FMcK2013b} R. Friedrich, J. McKay, {\it The $S$-transform in arbitrary dimensions}, arXiv (2013)
\bibitem{FMcK2015} R. Friedrich, J. McKay, {\it Homogeneous Lie Groups and Quantum Probability}, arXiv (2015)
\bibitem{F2017} R. Friedrich, {\it (Co)monads in Free Probability Theory}, arXiv (2017)
\bibitem{F2019} R.~Friedrich, {\it Generalised operations in free harmonic analysis}, Semigroup Forum, Springer, published online 31. January 2019
\bibitem{Fro} A. Fröhlich, {\it Formal Groups}, LNM, vol. {\bf 74}, Springer (1968)
\bibitem{GS}
J. Galambos, I. Simonelli, {\it Products of Random Variables, Applications to Problems of Physics and to Arithmetic Functions}, Pure and Applied Mathematics, Marcel Dekker, New York, (2004)
\bibitem{Haa1997}
U. Haagerup, {\it On Voiculescu's $R$- and $S$-Transforms for Free Non-Commuting Random Variables}, Free Probability Theory, Fields Institute Commun., Vol. {\bf 12}, AMS, 312pp, (1997)

\bibitem{Hz} M. Hazewinkel, {\it Witt vectors. Part 1}, revised version: 20 April 2008

\bibitem{Hz1978} M. Hazewinkel, {\it Formal Groups And Applications}, Academic Press, (1978)
  
\bibitem{H1956}
  F. Hirzebruch, {\it Neue topologische Methoden in der algebraischen
Geometrie}, Ergebnisse der Mathematik,  Springer Verlag, (1956)

\bibitem{Ho1968} T. Honda, {\it Formal Groups and Zeta-Functions}, Osaka J. Math. {\bf 5} (1968)
\bibitem{HO} A. Hora, N. Obata, {\it Quantum Probability and Spectral Analysis of Graphs}, Th. and Math. Physics, Springer, (2007)
\bibitem{I1979} L. Illusie, {\it Complexe de de Rham-Witt et cohomologie cristalline},
Annales scientifiques de l’É.N.S. $4^e$ série, tome 12, no {\bf 4}, p. 501-661 (1979).

\bibitem{KSU} T. Katsura, Y. Shimizu, K. Ueno, {\it Complex cobordism ring
and conformal field theory over $\mathbf{Z}$},  Math. Ann. 291,551-571 (1991)

\bibitem{KNTY} N. Kawamoto, Y. Namikawa, A. Tsuchiya,  Y. Yamada, {\it Geometric Realization of Conformal Field Theory on Riemann Surfaces}, Commun. Math. Phys. 116, 247–308 (1988).

\bibitem{K} S. V. Kerov, {\it Asymptotic Representation Theory of the Symmetric Group and Its Applications in Analysis}, Providence, Rhode Island: AMS (Translations of Mathematical Monographs, Vol. 219), (2003).

\bibitem{KY1988} A.A. Kirillov, D.V. Yuriev, {\it Representations of the Virasoro algebra by the orbit method}, JGP, Vol.~5, n.~{\bf 3}, (1988).

\bibitem{KLM} H. Kvinge, A.M. Licata, S. Mitchell,
{\it Khovanov’s Heisenberg category, moments in free probability, and shifted symmetric functions}, Algebraic Combinatorics, Centre Mersenne, Volume 2, issue {\bf 1} p. 49-74. doi : 10.5802/alco.32, (2019).\\
\url{https://alco.centre-mersenne.org/item/ALCO_2019__2_1_49_0/}

\bibitem{LTh} A. Lascoux, JY. Thibon, {\it Vertex operators and the class algebras of symmetric groups}, Journal of Mathematical Sciences 121: 2380. (2004) \\
\url{https://doi.org/10.1023/B:JOTH.0000024619.77778.3d}

\bibitem{Laz} M. Lazard, {\it Commutative Formal Groups}, LNM {\bf 443}, Springer Verlag, (1975)

\bibitem{Len} C. Lenart, {\it Formal Group-Theoretic Generalizations of the Necklace Algebra, Including a $q$-Deformation}, J. of Algebra. {\bf 199}, 703-732 (1998)
 

\bibitem{MSch2017} S. Manzel, M. Schürmann, {\it 
Non-commutative stochastic independence and cumulants}, Inf Dim Analysis, Quantum Probability and Related Topics Vol. 20, No. {\bf 02}, (2017) 

\bibitem{MN}
M. Mastnak, A. Nica, {\it Hopf Algebras and the Logarithm of the $S$-transform in Free Probability}, Trans. American Math. Soc., 
Vol. 362, Nr. {\bf 7}, p. 3705–3743  (2010)

\bibitem{McKS}
J. McKay,  A. Sebbar, {\it Replicable functions: An introduction}, Frontiers in Number Theory, Physics, and Geometry II: On Conformal Field Theories, Discrete Groups and Renormalization, Les Houches March 9-21, 2003, Springer, 1. ed., (2007)

\bibitem{MR} N. Metropolis, G-C. Rota, {\it Witt Vectors and the Algebra of Necklaces}, Adv. in Math. {\bf 50}, 95-125 (1983)

\bibitem{MJD} T. Miwa, M. Jimbo, E. Date, {\it Solitons: Differential Equations, Symmetries and Infinite Dimensional Algebras}, Cambridge University Press, (2000).

\bibitem{Mo}
J. Morava, {\it On the complex cobordism as a Fock representation},
LN in Mathematics, Vol. {\bf 1418}, Springer, (1990), 184-204.

\bibitem{SW1985} G. Segal, G. Wilson, 
{\it Loop groups and equations of KdV type},
Publications mathématiques de l'I.H.É.S., tome {\bf 61}, p. 5-65. (1985)

\bibitem{NS}
 A. Nica, R. Speicher, {\it Lectures on the Combinatorics of Free
Probability}, LMS LNS {\bf 335}, Cambridge University Press, (2006)


\bibitem{P} C. Pommerenke, {\it Univalent Functions}, Vandenhoeck $\&$ Ruprecht, Studia Mathematica / Mathematische Lehrbücher, Band XXV, Göttingen (1975)

\bibitem{S}
R. Speicher, {\it Free probability theory and non-crossing partitions},
S{\' e}minaire Lotharingien de Combinatoire, {\bf B39c}, 38pp., (1997)

\bibitem{Ta2001} L. Takhtajan,  {\it Free Bosons and Tau-Functions for Compact Riemann Surfaces and Closed Smooth Jordan Curves. Current Correlation Functions}, Lett. in Math. Phys. 56: 181–228, (2001)

\bibitem{Teo2003} 
L.-P. Teo, {\it Analytic Functions and Integrable Hierarchies-Characterization of Tau Functions}, Lett. Math. Phys. {\bf 64}: 75–92, (2003)

\bibitem{V1985}
  D.V. Voiculescu, {\it Symmetries of some reduced free product $C^*$-algebras, Operator Algebras and Their Connections with Topology and Ergodic Theory}, Lecture Notes in Mathematics, vol. {\bf 1132}, Springer Verlag, pp. 556-588, (1985)
  
\bibitem{V1986}
 D.V. Voiculescu, {\it Addition of Certain Non-commuting Random Variables},  J. Functional Analysis {\bf 66}, 323-346 (1986)

\bibitem{V1987}
  D.V. Voiculescu, {\it Multiplication of certain non-commuting random
variables}, J. Operator Theory, {\bf 18}, 223-235, (1987)

\bibitem{VDN}
  D.V. Voiculescu, K.J. Dykema, A. Nica {\it Free random variables: a noncommutative probability approach to free products with applications to random matrices, operator algebras, and harmonic analysis on free groups}, AMS, (1992)
  
\bibitem{V1998} D.V. Voiculescu, {\it The analogues of entropy and of Fisher's information measure in free probability theory; 
V. Noncommutative Hilbert Transforms}, Invent. math. {\bf 132}, (1998)

\bibitem{W1979} W. Waterhouse, {\it Introduction to Affine Group Schemes}, Springer (1979)

\bibitem{Wiki} 
Wikipedia,
\url{https://de.wikipedia.org/wiki/Plancksches_Strahlungsgesetz}


\end{thebibliography}
\end{document}